\documentclass[12pt]{article}
\usepackage[T1]{fontenc}
\usepackage{lmodern}
\usepackage[a4paper,hmargin=25mm,vmargin=23mm]{geometry}
\usepackage{amsmath,amssymb,amsthm,mathtools,booktabs,array,enumitem,microtype}
\usepackage{graphicx,placeins,needspace}
\usepackage[colorlinks=true,linkcolor=blue,citecolor=blue,urlcolor=blue]{hyperref}
\hypersetup{pdftitle={The Bernoulli problem for orthogonal partitions of a simplex},pdfauthor={Oleg R. Musin},pdfkeywords={Bernoulli quadrisection; orthogonal mass partition; simplex; tetrahedron; Dirichlet distribution; Borsuk--Ulam theorem; discriminant}}
\newtheorem{theorem}{Theorem}[section]
\newtheorem{lemma}[theorem]{Lemma}
\newtheorem{proposition}[theorem]{Proposition}
\newtheorem{corollary}[theorem]{Corollary}
\theoremstyle{definition}
\theoremstyle{remark}

\newcommand{\vol}{\operatorname{vol}}
\newcommand{\conv}{\operatorname{conv}}
\let\mathscr\mathcal
\setlist{itemsep=2pt,topsep=4pt}
\title{\Large The Bernoulli problem for orthogonal partitions of a simplex}
\author{Oleg R. Musin}
\date{}
\begin{document}
\maketitle

\begin{abstract}
Jacob Bernoulli's problem asks for two perpendicular lines dividing a
triangle into four equal areas. For an $n$-simplex, we study $n$ mutually
perpendicular hyperplanes and distinguish pairwise partitions, in which
every pair quarters the volume, from full partitions into $2^n$ equal parts.
We prove generic finiteness of full partitions in every dimension and of
pairwise partitions for tetrahedra. We construct simplices with exactly one
orthogonal partition in every dimension and exactly three in dimensions
$n\ge2$, all full. A regular $n$-simplex admits a full partition if the binary
expansion of $n+1$ has at most three ones. For tetrahedra, we classify
partitions through the centroid, determine which are full, and prove that
every nonconstant centered family is an isolated connected component of the
solution space. We also characterize partitions containing a face-parallel
or edge-containing cut. The proofs use affine configuration spaces,
Dirichlet projections, and axial stretching. Computer-assisted arguments
give a complete weighted-triangle count and an open set of tetrahedra
with at least seventeen pairwise partitions.
\end{abstract}



\section{Introduction}
\label{sec:introduction}

\subsection{The triangle problem and its solution map}
\label{subsec:triangle-introduction}

Jacob Bernoulli's problem asks for two perpendicular lines dividing a
triangle into four regions of equal area. In 1687 Bernoulli
\cite{Bernoulli} reduced the construction to an algebraic equation
of degree eight. Euler returned to the problem in two papers written
in 1779 and published in 1809 \cite{Euler729,Euler730}. Besides
existence, the problem has a natural counting question: how many
different pairs of lines can a given triangle admit?

For a triangle, a \emph{quadrisection} means an unordered pair of
unoriented perpendicular lines giving four equal areas. Symmetries
of the triangle may carry one quadrisection to another; these are
counted as different geometric configurations. The planar answer
is particularly simple.

\begin{theorem}[Triangle quadrisections; {\cite[Theorems 2 and 3]{Eberhart}}]
\label{intro:triangle-count}
Every nondegenerate triangle has one, two, or three quadrisections.
A scalene triangle has exactly one whose triangular cell lies on
its middle-length side, none on its longest side, and zero, one,
or two on its shortest side. The equilateral triangle has three.
\end{theorem}

The triangular cell in this description is bounded by the two cutting lines and a segment of the indicated side.

\begin{figure}[htbp]
 \centering
 \includegraphics[width=\linewidth]{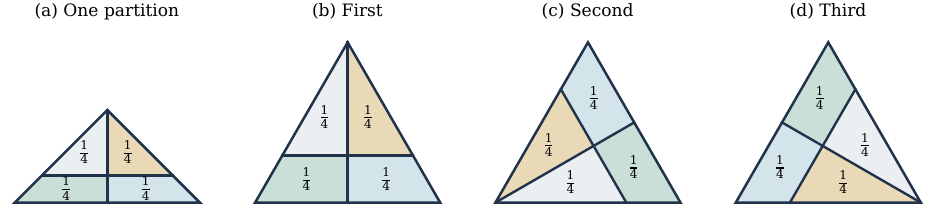}
 \caption{The unique quadrisection of a right isosceles triangle (a)
 and the three quadrisections of an equilateral triangle (b--d).
 Each label denotes a fraction of the area of its own triangle.
 These are separate pairs of perpendicular lines.}
 \label{fig:triangle-quadrisections}
\end{figure}

Figure~\ref{fig:triangle-quadrisections} displays the contrast between
a single configuration and three configurations of the same triangle.
For the isosceles triangle with vertices $(0,0),(1,0),(1/2,z)$, one
quadrisection is
\[
 x=\frac12,\qquad y=z\left(1-\frac1{\sqrt2}\right).
\]
The horizontal line cuts off a similar triangle of half the total
area, and the vertical symmetry axis bisects both parts. For
$z=1/2$ this is the only quadrisection. For $z=\sqrt3/2$, rotating
both lines through $120^\circ$ and $240^\circ$ about the centroid
gives the other two. Their three intersection points are distinct.

\begin{figure}[htbp]
 \centering
 \includegraphics[width=\linewidth]{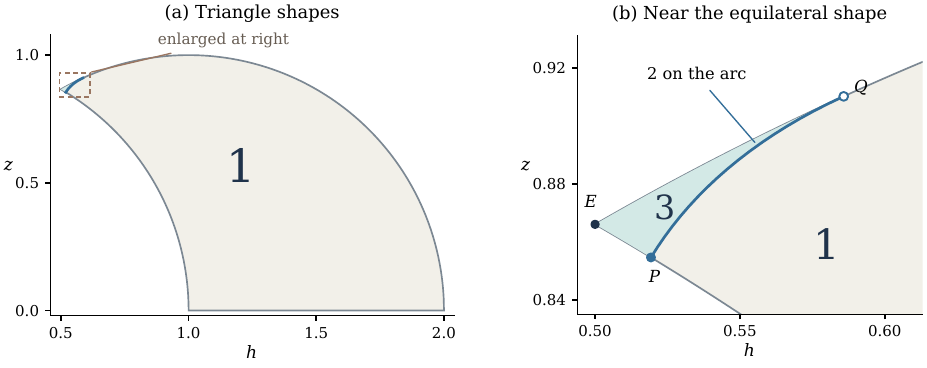}
 \caption{The number of quadrisections in the triangle shape domain
 \eqref{eq:intro:triangle-domain}; the right panel enlarges the region
 near $E$. On either isosceles boundary the count is three from $E$ up to the separating arc and one beyond it. The separating arc has two solutions, including
 $P=(175/337,288/337)$, with the sole exception
 $Q=(2-\sqrt2,\sqrt{2\sqrt2-2})$, which has one. The equilateral
 shape $E$ has three. This is the uniform-area triangle problem.}
 \label{fig:triangle-count-map}
\end{figure}

The count can also be displayed over the space of triangle shapes.
Normalize a side of middle length to the segment from $(0,0)$ to
$(1,0)$ and put the remaining vertex at $(h,z)$, with $h\ge1/2$.
The resulting shape domain is
\begin{equation}\label{eq:intro:triangle-domain}
 \Omega=\{(h,z):h\ge1/2,\ z>0,\ h^2+z^2\ge1,
                         \ (h-1)^2+z^2\le1\}.
\end{equation}
Its curved boundary consists of isosceles shapes; the two arcs meet
at the equilateral shape $E=(1/2,\sqrt3/2)$. The map in
Figure~\ref{fig:triangle-count-map} is redrawn from the area and
perpendicularity equations of \cite[Section 4.2, Theorem 3]{Eberhart}. It shows a region
with three solutions next to $E$ and a region with one solution.
Along their separating arc the count is two, except at the indicated
endpoint $Q$, where it is one. Thus even in dimension two, the
existence theorem and the description of all solutions are different
parts of the problem.

\subsection{Two problems for a simplex}
\label{subsec:two-problems}

Bernoulli's two perpendicular lines give four equal areas. Taking one
perpendicular cut for each dimension leads directly to the following question.

\paragraph{Problem 1: full orthogonal equipartitions.}
For which $n$-dimensional simplices do $n$ mutually perpendicular
hyperplanes divide the simplex into $2^n$ equal-volume parts?
Describe such simplices and their full partitions. In particular,
determine which regular simplices admit a full partition and how
the property changes under deformation of the simplex.

A related existence theorem provides the starting point for our second question.
The result was stated by Makeev~\cite{Makeev}; a complete proof is given
in~\cite[Theorem~1.2]{MusinOrth}.

\Needspace{10\baselineskip}
\begin{theorem}[Orthogonal four-partitions]
\label{intro:orthogonal-existence}
Let $n\ge2$, and let $\mu$ be a finite nonzero Borel measure on
$\mathbb R^n$ assigning measure zero to every affine hyperplane.
There exist $n$ mutually perpendicular affine hyperplanes
$H_1,\ldots,H_n$ such that, for every $i\ne j$ and every choice
$\sigma,\tau\in\{+,-\}$,
\[
 \mu(H_i^\sigma\cap H_j^\tau)=\frac14\mu(\mathbb R^n),
\]
where $H_i^+$ and $H_i^-$ are the two half-spaces bounded by $H_i$.
\end{theorem}

For simplex volume, call the resulting system an \emph{orthogonal partition}:
each plane bisects the volume and every pair quarters it. In dimension one
only the median condition is needed. A partition is \emph{full} if each of the
$2^n$ orthants has the same volume. In a triangle the two conditions agree;
for a tetrahedron fullness imposes one further equation, as in~\eqref{eq:gen:triple-moment}.

\paragraph{Problem 2: the number of orthogonal partitions.}
How many orthogonal partitions does a given simplex have? Describe
the finite counts, the shapes carrying continuous families, and
the changes in the solution set as the shape varies. Existence
in this problem is supplied by Theorem~\ref{intro:orthogonal-existence};
the remaining issue is to understand all configurations.

We count unordered collections of unoriented planes. A symmetry of the
simplex may carry one partition to another, as in the equilateral example.

\subsection{Main results and open questions}
For Problem 1 we obtain a sufficient condition in every dimension and an
exact criterion for partitions through the centroid in dimension three.

\begin{theorem}\label{intro:regular}
If the binary expansion of $n+1$ contains at most three ones, the regular
$n$-simplex admits a full orthogonal partition.
\end{theorem}

The binary condition is sufficient. Whether it is necessary is open:
the first case it leaves undecided is the regular $14$-simplex, with
$15=8+4+2+1$ vertices. For the regular tetrahedron the theorem gives existence, since $3+1=2^2$.
In dimension three we can say more. We call a partition \emph{centered}
if all its planes pass through the centroid. The next theorem characterizes
the tetrahedra admitting a centered full partition and counts these partitions.

\begin{theorem}\label{intro:center}
A tetrahedron admits a full orthogonal partition through its centroid
if and only if four edges forming a cycle have the same length.
Among such partitions there is exactly one for each nonregular
tetrahedron satisfying this condition and exactly three for the regular tetrahedron.
\end{theorem}

The criterion concerns partitions through one specified point. Even for the
regular tetrahedron it is open whether the three centered octasections exhaust
all full partitions.

For Problem 2 the following counts include every orthogonal partition
of the given simplex.

\begin{theorem}\label{intro:counts}
For each $n\ge1$ there is an $n$-simplex with exactly one orthogonal
partition. For each $n\ge2$ there is an $n$-simplex with exactly
three. In both constructions every partition is full.
\end{theorem}

The existence proof in Section~\ref{sec:general} works in all dimensions.
In dimension three, Theorem~\ref{thm:weighted-atlas-certified} also identifies
an explicit planar shape used in the construction.
These are exact counts on particular simplices. Which other finite counts occur
is open already for tetrahedra; Section~\ref{subsec:seventeen-computations}
exhibits an open set with at least seventeen solutions.

\begin{figure}[htbp]
\centering
\includegraphics[width=.95\linewidth]{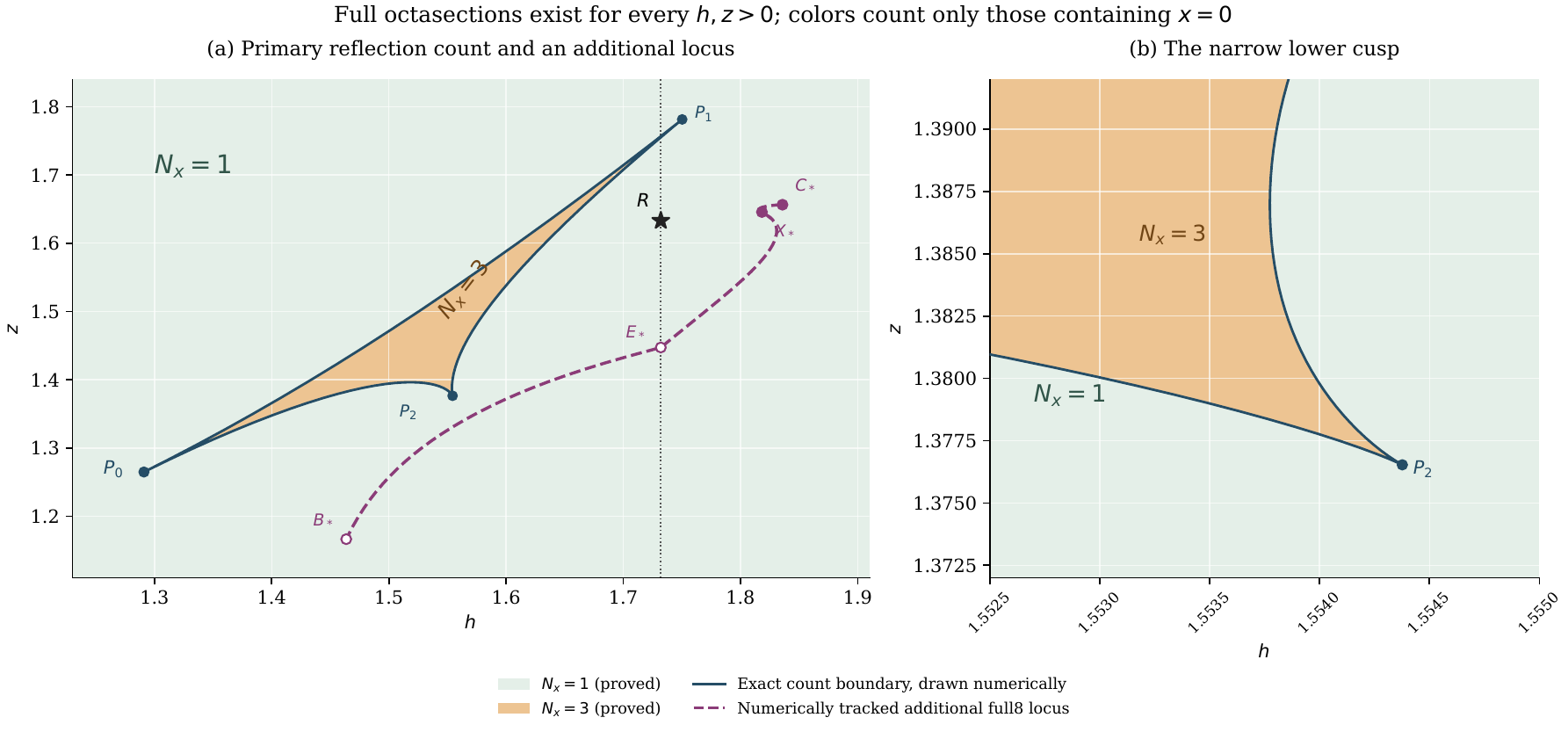}
\caption{Full octasections containing $x=0$. The solid curve separates
the regions $N_x=1$ and $N_x=3$. On the curve $N_x=2$ except at its three
tips, where $N_x=1$. The right panel enlarges the lower cusp.
The dashed curve tracks an additional branch numerically and does not enter
the color count. The regular tetrahedron is marked $R$.}
\label{fig:hz}
\end{figure}

Figure~\ref{fig:hz} illustrates how counts change in a two-parameter family.
Take vertices $(\pm1,0,0),(0,h,0),(0,(h^2-1)/(2h),z)$ with $h,z>0$.
To follow a single construction, count only the partitions containing the
mirror $x=0$. By Theorem~\ref{thm:c11-edge} these are full and correspond
to quadrisections of the projected triangle. Scaling its coordinates by $1/h$
gives $a=(h^2-1)/(2h^2)$ and $K=a^2+(z/h)^2$ in Section~\ref{subsec:weighted-map}.
In higher dimensions the same planar construction gives regions with one or
three solutions for the examples of Section~\ref{sec:general}; the weights
vary with the dimension, and a common global shape diagram has not been proved.

\begin{theorem}\label{intro:finite}
For a generic Euclidean $n$-simplex the set of full orthogonal
partitions is finite, possibly empty. For a generic tetrahedron
the set of pairwise orthogonal partitions is finite and its cardinality
is locally constant. Here generic means outside a proper real algebraic
subset of the positive definite metrics on a fixed affine simplex.
\end{theorem}

For full partitions even the generic existence question remains open in
dimension three. For pairwise partitions, generic finiteness in dimensions
$n\ge4$ is also open. The exceptional case below shows why a genericity
hypothesis matters.

\begin{theorem}\label{intro:circles}
The centered orthogonal partitions of a regular tetrahedron form three circles
with one common point. Each circle contains exactly one full partition.
\end{theorem}

For an arbitrary tetrahedron, Theorems~\ref{c11:centered-pairwise}--\ref{thm:v33:centered-positive-components}
give the exact criterion for a nonconstant centered family: the three bimedians
must be mutually perpendicular and at least two must have equal length.
A bimedian joins the midpoints of opposite edges. For a nonregular tetrahedron
satisfying the criterion the family is one circle. In both cases it is an
isolated component. The full partitions away from the centroid remain to be
classified. Proposition~\ref{prop:liftfamily} extends the existence of continuous
families to regular simplices in all higher dimensions. Whether a fixed tetrahedron has a nonconstant family outside the centered locus remains open.

\subsection{Median frames and the volume equations}
\label{subsec:balanced-cuts}

The input for a simplex is the finite set of vertices
$\mathcal V=\{v_0,\ldots,v_n\}$. For any proposed cuts, each cell is a
polytope determined by these vertices and the hyperplanes.
Its volume can be computed from a finite list of intersection
points: triangulate the cell and add the determinant formulas for
the resulting simplex volumes. Thus the equations used below
start from finitely many coordinates. The balanced quantities
are cell volumes, rather than the number of original vertices in
each cell.

Fix a normal direction and move a perpendicular hyperplane through
the simplex. There is a unique position giving equal volume on
the two sides. An orthonormal frame therefore determines all
median hyperplanes, leaving one equation per pair of normals.
For a tetrahedron there are three such equations. Once they hold,
one additional equation decides whether the eight volumes are equal.

The classical Borsuk--Ulam theorem \cite{Borsuk} forces a zero of
a continuous odd map from the $d$-sphere to $\mathbb R^d$.
Its manifold and frame versions \cite{Musin2012,MusinStiefel,MusinOrth}
apply to the map of pairwise volume imbalances on median frames.
Reversing or permuting normals acts on these coordinates; the corresponding
equivariant zero theorem gives Theorem~\ref{intro:orthogonal-existence}.

The continuation method of~\cite{MusinNumerical} follows known solutions
as the simplex is deformed. Exact volume formulas and rigorous error bounds verify the resulting candidates.

Section~\ref{sec:general} proves generic finiteness results and constructs simplices
with exactly one or three partitions. Section~\ref{sec:tetrahedra} classifies the
centered partitions of a tetrahedron and determines which give eight equal volumes.
Section~\ref{sec:computations} counts quadrisections of a weighted triangle as its
shape varies, exhibits an open family of tetrahedra with at least seventeen
orthogonal partitions, and concludes with open questions.

\FloatBarrier

\Needspace{7\baselineskip}
\section{General results in arbitrary dimension}
\label{sec:general}

We prove Theorems~\ref{intro:regular}, \ref{intro:counts}, and
\ref{intro:finite}. Two arguments are involved. For finiteness, first
allow arbitrary angles between the cuts and bound the dimension of the
resulting space. Orthogonality then gives linear equations on the
metric. For pairwise partitions of a tetrahedron the decisive fact is
that two quartering planes determine the third affine plane uniquely.

For existence and exact counts we start with a weighted triangle.
Grouping the barycentric coordinates of a simplex into three blocks
projects its volume onto the triangle. Each planar quartering lifts
to a full partition. For the regular simplex, blocks of power-of-two
size give perpendicular cuts within each block, proving
Theorem~\ref{intro:regular}. For the exact counts a further step is
essential: stretching one axis at a time forces every partition to
contain the prescribed cuts. The count therefore agrees with the
planar count.

\subsection{Configuration spaces and the finiteness statements}
\label{gen:definitions}

Let $\mu$ be a probability measure on a full-dimensional convex body in
$\mathbb R^n$, with positive interior density, and let $X$ have law $\mu$.
A \emph{median hyperplane}
has mass $1/2$ on each side. For a nonzero normal $a$, write $h_\mu(a)$
for the unique median of $a\cdot X$; thus its median plane is
$a\cdot x=h_\mu(a)$. For uniform volume on $P$ we write $h_P$.
An orthogonal partition consists of $n$
mutually perpendicular median hyperplanes $L_i=0$ such that every pair
quarters the mass. It is \emph{full} if every one of its $2^n$ sign cells
has mass $2^{-n}$. Medianity is required also for $n=1$: in that case
the partition is the unique median point. We disregard the order of
the planes and the signs of their normals, but not symmetries of the
underlying body. Denote these two configuration sets by
$\mathcal P_n(\mu)$ and $\mathcal F_n(\mu)$.

Let $p_\varepsilon$ be the mass of the sign cell
$\{\varepsilon_iL_i>0\text{ for all }i\}$, where $\varepsilon_i\in\{-1,1\}$.
For oriented defining forms put $s_i=\operatorname{sgn}L_i$,
$F_I=\mathbb E_\mu\prod_{i\in I}s_i$, and $F_\varnothing=1$.
Expanding $\prod_i(1+\varepsilon_i s_i)/2$ gives the sign-moment identity
\begin{equation}\label{eq:gen:walsh}
 p_\varepsilon=2^{-n}\sum_{I\subseteq\{1,\ldots,n\}}
       F_I\prod_{i\in I}\varepsilon_i.
\end{equation}
Thus the pairwise conditions are $F_I=0$ for $1\le|I|\le2$,
whereas fullness requires every nonempty moment to vanish.
We abbreviate $F_{\{i,j\}}$ to $F_{ij}$. The \emph{reduced pairwise moment map}
assigns to an orthonormal frame its $n(n-1)/2$ pair moments, after placing
each cut at its median. A partition is \emph{regular} if the Jacobian of
this map in local frame coordinates is invertible. This condition is
independent of the chosen coordinates and of the order and orientation
of the cuts. In dimension two we also call a regular quartering \emph{simple}:
the derivative of its single moment with respect to the frame angle is nonzero.
For median cuts, \eqref{eq:gen:walsh} gives the two identities used repeatedly below:
\begin{equation}\label{eq:gen:pair-moment}
 F_{ij}=4\Pr(L_i>0,L_j>0)-1,
\end{equation}
and, for a pairwise triple,
\begin{equation}\label{eq:gen:triple-moment}
 p_\varepsilon=\frac{1+\tau\varepsilon_1\varepsilon_2\varepsilon_3}{8},
 \qquad \tau=F_{\{1,2,3\}}=8p_{+++}-1.
\end{equation}
In particular, a pairwise triple is full exactly when $\tau=0$.

For the affine results, $P$ is a full-dimensional convex polytope with
uniform volume. Write $\operatorname{Sym}_n^+$ for the cone of symmetric
positive definite $n\times n$ matrices. An inverse metric $Q$ acts on
normal covectors: orthogonality means $a_i^{\mathsf T}Qa_j=0$.
When such a metric is specified, we write
$\mathcal P_n(P,Q)$ and $\mathcal F_n(P,Q)$ for the corresponding
orthogonal configuration sets. Let $\mathcal A_k(P)$ be the space of ordered,
unoriented affine $k$-tuples with $2^k$ equal cells, and let
$\mathcal B_3(P)$ be the corresponding space of pairwise quartering
triples. No orthogonality is imposed in these two spaces.
All dimensions below are semialgebraic dimensions.

Equal-volume conditions determine the affine arrangements; the metric
determines which are orthogonal. The following dimension estimate does
not assume independence of the derivatives of the volume equations.

\begin{theorem}[Affine dimension and generic full finiteness]
\label{gen:affine-dimension}\label{gen:generic-full}
For $0\le k\le n$,
\begin{equation}\label{eq:gen:affine-dimension}
 \dim\mathcal A_k(P)\le kn-\frac{k(k+1)}2.
\end{equation}
Fixing the affine model and varying its inverse Euclidean metric
$Q\in\operatorname{Sym}_n^+$, a cone of dimension $m=n(n+1)/2$,
the locus whose full orthogonal solution space has dimension at
least $d\ge1$ has dimension at most $m-d$. In particular, outside
a proper real algebraic subset, every full orthogonal solution set
is finite, possibly empty. For a simplex this is generic finiteness
with respect to its Euclidean shape.
\end{theorem}

\begin{theorem}[Quartering partners and completion]
\label{gen:projective-fiber}\label{gen:completion3}
For $n\ge2$, the quartering partners of a fixed median $L=0$
are parametrized by the projective restriction of their normals to
$\ker dL$. This is a $C^1$ semialgebraic diffeomorphism with
$\mathbb {RP}^{n-2}$. Consequently $\mathcal A_2(P)$ is a compact
$C^1$ semialgebraic manifold of dimension $2n-3$, diffeomorphic to
the space of ordered perpendicular unoriented directions.

When $n=3$, every quartering pair has a unique affine third plane
completing it to a pairwise partition. Forgetting that plane gives
a $C^1$ semialgebraic diffeomorphism
$\mathcal B_3(P)\cong\mathcal A_2(P)$. The three median and three
pairwise equations have rank six throughout $\mathcal B_3(P)$.
These assertions include planes through vertices or edges.
\end{theorem}

\begin{theorem}[Generic pairwise finiteness in dimension three]
\label{gen:generic-pairwise}
For $P\subset\mathbb R^3$ and $1\le d\le3$,
\begin{equation}\label{eq:gen:pairwise-exceptional}
 \dim\{Q\in\operatorname{Sym}_3^+:
                  \dim\mathcal P_3(P,Q)\ge d\}\le6-d.
\end{equation}
Outside a proper real algebraic subset, all pairwise orthogonal
partitions are finite in number and regular, and their count is
locally constant. Regularity means that the three median, three
pairwise, and three orthogonality equations have rank nine in
affine-plane coordinates.
\end{theorem}

\subsection{Dirichlet measures and explicit full partitions}
\label{gen:dirichlet}

The Dirichlet law $\operatorname{Dir}(\alpha_1,\ldots,\alpha_r)$ has density
proportional to $\prod_i\lambda_i^{\alpha_i-1}$ on
$\{\lambda_i\ge0,\ \sum_i\lambda_i=1\}$, where $\alpha_i>0$.
Its two-coordinate case gives the beta law $\operatorname{Beta}(a,b)$,
with density proportional to $t^{a-1}(1-t)^{b-1}$ on $0<t<1$.

For $n\ge2$, a uniform point of an $n$-simplex has barycentric law
$\operatorname{Dir}(1,\ldots,1)$. Partition its vertex indices into
three blocks $I_1,I_2,I_3$ of positive sizes $p_1,p_2,p_3$, and choose
a rooted binary tree in each block. At an internal node $v$, let
the two child sets $I_v^-,I_v^+$ have sizes $a_v,b_v$, respectively, and let $\theta_v$ be the
median of $\operatorname{Beta}(a_v,b_v)$. Put
\[
 L_v(\lambda)=(1-\theta_v)\sum_{j\in I_v^-}\lambda_j
                    -\theta_v\sum_{j\in I_v^+}\lambda_j.
\]
There are $n-2$ such forms. Choose positive lengths $\ell_v$ and a
nondegenerate triangle $\Delta=\conv(B_1,B_2,B_3)$, and set
\begin{equation}\label{eq:gen:lift-vertices}
 P_j=((\ell_vL_v(e_j))_v,B_i)\in\mathbb R^{n-2}\oplus\mathbb R^2
 \quad(j\in I_i),\qquad T=\conv(P_1,\ldots,P_{n+1}).
\end{equation}

\begin{theorem}[The planar lift]\label{gen:planar-lift}
The vertices in \eqref{eq:gen:lift-vertices} are affinely independent.
Full orthogonal partitions containing the coordinate cuts $x_v=0$
are in bijection with perpendicular quarterings of $\Delta$ for
the probability density
\begin{equation}\label{eq:gen:dirichlet-density}
 \rho(Y)=\frac{n!}{2\operatorname{area}(\Delta)\prod_i(p_i-1)!}
                       \prod_{i=1}^3w_i(Y)^{p_i-1},
\end{equation}
where $w_i$ are barycentric coordinates. Such quarterings always
exist, and
\[
 \vol_n(T)=\frac{2\operatorname{area}(\Delta)}{n!}\prod_v\ell_v.
\]
For fixed trees and block sizes these simplices form an
$n$-dimensional family up to similarity.
\end{theorem}

\begin{corollary}\label{gen:regular-series}
The regular $n$-simplex has a full orthogonal partition whenever
the binary expansion of $n+1$ has at most three nonzero digits.
\end{corollary}

\subsection{Global rigidity under axial stretching}
\label{gen:stretching}

Let $m\ge1$, and let $K\subset\mathbb R\times\mathbb R^m$ be a bounded
full-dimensional convex polytope with a polynomial probability
density $\rho$, positive inside. Assume that $X=0$ is conditionally
median on every interior projection fiber:
\begin{equation}\label{eq:gen:conditional-fairness}
 \int_{X>0}\rho(X,Y)\,dX=\int_{X<0}\rho(X,Y)\,dX.
\end{equation}
Let $\mu$ be the projected measure and let $\nu_L$ be the image
of $\rho$ under $(X,Y)\mapsto(LX,Y)$.

Condition~\eqref{eq:gen:conditional-fairness} says that the first
sign is fair even after the projected point $Y$ is known. Equivalently,
the conditional distribution of $Y$ on either side of $X=0$ is the
same measure $\mu$. The two halves of each fiber need not be mirror images. Stretching changes
their lengths while preserving this conditional probability identity.
Agreement of \emph{local solution germs} below means that the correspondence
identifies the solution sets in neighborhoods of corresponding partitions.

\begin{theorem}[Global axial rigidity]\label{gen:global-stretch}
For all sufficiently large $L$, every $(m+1)$-plane orthogonal
partition of $\nu_L$ contains the physical plane $x=0$. Deleting
that plane gives bijections
\[
 \mathcal P_{m+1}(\nu_L)\cong\mathcal P_m(\mu),\qquad
 \mathcal F_{m+1}(\nu_L)\cong\mathcal F_m(\mu).
\]
The inverse adds $x=0$. Corresponding local solution germs agree,
and the rank of either reduced moment map increases by exactly $m$.
No finiteness or regularity assumption on the base solution set is
required.
\end{theorem}

\begin{corollary}[Global transfer by integer splitting]
\label{gen:integer-split}
Let a simplex have positive integer Dirichlet weights, and let
a vertex $q$ have weight $r+s$, $r,s\ge1$. If $\theta$ is the
median of $\operatorname{Beta}(s,r)$, replace it by
\[
 q-a\theta e\quad\hbox{of weight }r,\qquad
 q+a(1-\theta)e\quad\hbox{of weight }s,
\]
where $e$ is a new perpendicular direction. For every sufficiently
large $a$, this preserves the entire pairwise and full solution
sets, their local germs, and the rank increment of
Theorem~\ref{gen:global-stretch}. Iterating until all weights are
one gives a uniform simplex of dimension $N-1$, where $N$ is the
original total weight.
\end{corollary}
\begin{proof}
If $\lambda$ is the aggregate coordinate, Dirichlet aggregation
gives a ratio $t\sim\operatorname{Beta}(s,r)$ independent of the
projected point, and the new coordinate is $a\lambda(t-\theta)$.
Its sign is conditionally fair and its projection is the original
weighted simplex. Integer weights give polynomial density positive
inside, so the theorem applies. Fix each new length after all
previous geometry has been fixed; the total weight remains $N$.
\end{proof}

\subsection{Simplices with exactly one or three partitions}
\label{gen:exact}

\begin{theorem}\label{gen:exact-counts}
For every $n\ge1$ there is a Euclidean $n$-simplex with exactly
one orthogonal pairwise partition, which is full. For every $n\ge2$
there is one with exactly three, all full. At every constructed
configuration the reduced pairwise differential is invertible and
the reduced full moment differential is injective.
For $n=1$ these are statements about the zero-dimensional frame space.
\end{theorem}

Theorem~\ref{intro:finite} follows from Theorems~\ref{gen:generic-full}
and~\ref{gen:generic-pairwise}. Corollary~\ref{gen:regular-series}
gives Theorem~\ref{intro:regular}. Theorem~\ref{gen:exact-counts}
gives Theorem~\ref{intro:counts}.

\subsection{Proofs of the affine structure results}
\label{gen:affine-proofs}

We first record the geometric and differential facts shared by the proofs.

\begin{lemma}[Sign cells and a common intersection]\label{gen:sign-cells}
Let $L_1,\ldots,L_k$ be affine forms on a full-dimensional convex body $K$.
All $2^k$ open sign cells meet $\operatorname{int}K$ if and only if the
normals are independent and the common zero flat meets $\operatorname{int}K$.
\end{lemma}
\begin{proof}
Consider the convex image $C=L(\operatorname{int}K)\subset\mathbb R^k$.
If every orthant meets $C$, then $0\in\operatorname{int}C$: otherwise
separation by a hyperplane through zero would exclude an orthant.
This gives full rank and an interior common zero. Conversely, at such a
zero the full-rank affine map sends an interior neighborhood onto a
neighborhood of zero, meeting every orthant.
\end{proof}

\begin{lemma}[Ordering of bisectors]\label{gen:ordering}
Let a convex body carry a density positive almost everywhere in its interior.
Two affine forms that bisect its mass cannot satisfy $f>g$ throughout the interior.
\end{lemma}
\begin{proof}
The positive half of $g$ lies in the positive half of $f$.
At an interior zero of $g$ we have $f>0$, so their difference contains
an open set on the negative side of $g$ and has positive mass.
\end{proof}

\begin{lemma}[Section derivatives and median reduction]\label{gen:cut-flux}
Let $d\mu=\rho(x)\,dx$ on a convex polytope, with $\rho$ continuous,
bounded and positive in the interior. For an interior cut $L=a\cdot x-h$
and a smooth variation of its coefficients,
\begin{equation}\label{eq:gen:cut-flux}
 \frac{d}{dt}\mu(L>0)
 =\frac1{\|a\|}\int_{L=0}(\dot a\cdot x-\dot h)\rho(x)\,dS.
\end{equation}
For independent moving half-spaces, sum this flux over each moving
section restricted by the others. These derivatives are continuous.
Median offsets are $C^1$, with
\begin{equation}\label{eq:gen:median-derivative}
 Dh_\mu(a)[v]
 =\frac{\int_{a\cdot x=h_\mu(a)}(v\cdot x)\rho(x)\,dS}
        {\int_{a\cdot x=h_\mu(a)}\rho(x)\,dS}.
\end{equation}
Regularity of the reduced pairwise map is equivalent to full rank $n^2$
of the median, pairwise and orthogonality equations in affine-plane coordinates.
\end{lemma}
\begin{proof}
Integrate along the normal to a moving plane. Dividing the swept mass by
the displacement gives the displayed section integral by dominated convergence.
Intersections of two moving slabs have volume $O(t^2)$ by independence.
On each interior section the other boundaries have section measure zero;
dominated convergence in a local graph chart gives continuity even at
vertex incidences. The offset derivative is strictly negative, so the
median formula follows by implicit differentiation. At an orthogonal frame,
the $n$ median constraints eliminate the offsets, and the $n(n-1)/2$
independent orthogonality constraints leave precisely the frame variables.
The remaining $n(n-1)/2$ pairwise equations are the reduced map, proving
the rank equivalence.
\end{proof}

\begin{lemma}[Continuation of regular solutions]\label{gen:continuation}
Distinct regular roots of a square $C^1$ system persist on disjoint local
branches under small parameter changes. If all roots above a parameter
are regular and the solution projection is proper nearby, their total
number is locally constant.
\end{lemma}
\begin{proof}
The implicit function theorem gives the branches. A compact discrete fiber
is finite; any sequence of additional roots approaching its parameter
would, by properness, accumulate on one of those branches.
\end{proof}

\paragraph{Restriction of the next normal.}
The relevant spaces and volume functions are semialgebraic. Indeed,
choose independent active facet equations for every clipped vertex,
retain their determinant and feasibility sign conditions, and
triangulate each clipped cell in a fixed vertex order. On each of
the resulting finitely many strata, vertices and volumes are rational
functions of the plane coefficients. Singular incidences are retained
as separate strata, and lower-dimensional cells have zero volume.

In a full $k$-tuple the normals are independent by
Lemma~\ref{gen:sign-cells}. Fix its forms $L_1,\ldots,L_k$, and set
$W=\bigcap_i\ker dL_i$. A simultaneous bisector $M=0$ of all its
cells has $dM|_W\ne0$. Otherwise $M=c+\sum a_iL_i$ is strictly
one-signed on a cell obtained by choosing the signs of its nonzero
terms to agree with $c$ (either choice when $c=0$).

Moreover, the projective restriction determines $M=0$ uniquely.
Normalize two candidates to have the same restriction. Their
difference is $c+\sum a_iL_i$ and, unless zero, has a strict fixed
sign on some open cell. Lemma~\ref{gen:ordering} applied to that cell excludes this possibility. Hence the fiber of
$\mathcal A_{k+1}\to\mathcal A_k$ injects semialgebraically into
$\mathbb P(W^*)$, of dimension $n-k-1$. The semialgebraic
fiber-dimension theorem \cite{Hardt,CosteRealSets} and summation give
\eqref{eq:gen:affine-dimension}.

If the normal rows of an affine full partition form $A$, all inverse
metrics making them orthogonal are precisely
\begin{equation}\label{eq:gen:metric-fiber}
 Q=A^{-1}\operatorname{diag}(d_1,\ldots,d_n)A^{-\mathsf T},
 \qquad d_i>0.
\end{equation}
The incidence space of full partitions and compatible metrics
therefore has dimension at most $n(n-1)/2+n=m$. Its projection to
the metric cone gives the exceptional-dimension bound. An infinite
semialgebraic set has positive dimension; the real Zariski closure
of a semialgebraic set has the same dimension
\cite[Theorem~2.2]{CosteRealSets}. This proves
Theorem~\ref{gen:generic-full}.

For a single previous cut, every projected direction actually occurs, as we now show.

\paragraph{The quartering fiber.}
Write $x=L$, and choose a transverse form $t$ with a prescribed
oriented restriction to $\ker dL$. The candidate partner has the
unique normalized form $H=t+ax+c$. The median $c(a)$ is unique and $C^1$ by Lemma~\ref{gen:cut-flux}.
Its section $D$ meets both open sides of $x=0$: otherwise one half-body
lies on one side of $H=0$ and exhausts its half-volume, although that same
side also meets the other half-body in an open set.

Let $\bar x_D$ be the section mean and put
$C(a)=\mathbb E[\operatorname{sgn}x\operatorname{sgn}H]$. Equations~\eqref{eq:gen:cut-flux}--\eqref{eq:gen:median-derivative} give
\[
 c'(a)=-\bar x_D,\qquad
 C'(a)=\frac{2\operatorname{area}_{n-1}(D)}
              {\vol(P)\|\nabla H\|}
       \operatorname{Cov}_D(\operatorname{sgn}x,x)>0.
\]
The conditional means of $x$ on its two sign halves are strictly ordered.
Since $ax$ has median zero, $|c(a)|\le\sup_P|t|$. Hence $C(a)\to\pm1$
as $a\to\pm\infty$ by dominated convergence. Its unique zero depends
$C^1$ on the transverse direction and is semialgebraic by its graph.

An auxiliary inner product now identifies $\mathcal A_2(P)$ with
the flag space
\[
 \{(\ell_1,\ell_2)\in(\mathbb {RP}^{n-1})^2:
                                      \ell_1\perp\ell_2\}:
\]
the first direction specifies its median, and the second specifies
the projected normal of its partner. The base and fiber dimensions
are $n-1$ and $n-2$. This proves the general part of
Theorem~\ref{gen:projective-fiber}.

\begin{lemma}[A section matrix]\label{gen:section-matrix}
For a convex polygon $D$ with $0\in\operatorname{int}D$, set
\begin{equation}\label{eq:gen:section-matrix}
 A_D=\int_D
 \begin{pmatrix}1\\\operatorname{sgn}x\\\operatorname{sgn}y\end{pmatrix}
 (1,x,y)\,dx\,dy.
\end{equation}
All principal minors of $A_D$ are positive.
\end{lemma}
\begin{proof}
The diagonal entries
are positive. The minors involving the first row are strict
covariances, as above. The remaining two-by-two minor is positive
because $|\int_D\operatorname{sgn}x\,y|<\int_D|y|$ and its
$x,y$ counterpart; all four quadrants meet the interior.

Here is a proof of the full determinant, which is the essential
point. Let $M,\mu$ be the area and centroid of $D$, and write
$u_\theta=(\cos\theta,\sin\theta)$ and
$t_\theta=(-\sin\theta,\cos\theta)$. Put
$w(\theta)=\int_D\operatorname{sgn}(u_\theta\cdot X)(X-\mu)\,dX$.
If $M_\pm,c_\pm$ are the two half-polygon areas and centroids,
\[
 w=\frac{2M_+M_-}{M}(c_+-c_-),\qquad w\cdot u_\theta>0.
\]
Write their common chord as $[-a t_\theta,b t_\theta]$.
The centroid line meets it at $r_0t_\theta$, where
\[
 (b-2a)/3\le r_0\le(2b-a)/3.
\]
Indeed, make that intersection the origin and the centroid line
horizontal. A polygon whose extreme vertical coordinates are
$y_-,y_+$ has centroid between $(2y_-+y_+)/3$ and
$(y_-+2y_+)/3$, by triangulating from an extreme vertex.
Convexity prevents both half-polygons from extending below the
lower chord endpoint, or both from extending above the upper one.
The shifted chord lengths are $a+r_0$ and $b-r_0$; the centroid
bounds give $b-r_0\le2(a+r_0)$ and $a+r_0\le2(b-r_0)$,
which are the displayed middle-third inequalities.

Endpoint differentiation of the rotating chord gives
\[
 w'=2(A t_\theta-B\mu),\quad
 A=(a^3+b^3)/3,\quad B=(b^2-a^2)/2.
\]
The formula holds at vertex crossings as well: integrate first in
polar radius and differentiate the two angular endpoints. Since
$\mu$ lies on the centroid line,
\[
 \det(w,w')=2(w\cdot u_\theta)(A-Br_0)>0,
 \qquad A-Br_0\ge\frac{\min(a,b)(a+b)^2}{6}.
\]
The argument of $w(\theta)$ is strictly increasing, lies within
$\pi/2$ of $\theta$, and increases by $\pi$ over a half-turn.
Thus $\det(w(0),w(\pi/2))>0$. Column elimination in
\eqref{eq:gen:section-matrix} yields
$\det A_D=M\det(w(0),w(\pi/2))>0$.
\end{proof}

\paragraph{Pairwise triples in three dimensions.}
For a pairwise triple, \eqref{eq:gen:triple-moment} has $|\tau|<1$:
a small ball about an interior point of one median plane, away from the
other two, meets positive-volume cells of both parities. All eight cells
therefore have positive volume. Lemma~\ref{gen:sign-cells} gives independent
normals and an interior common point.

For uniqueness, choose coordinates with the given pair $x=0,y=0$.
Its common axis crosses the interior; let its chord endpoints be
$z_-<z_+$. The ham-sandwich theorem, applied to volume on $P$,
$P\cap\{x>0\}$, and $P\cap\{y>0\}$, gives a completing plane.
Every completer has a unique form
$z+c+ax+by=0$, with $(c,a,b)$ in the rectangle
$\Omega=(-z_+,-z_-)\times\mathbb R^2$. Define
\[
 \Psi(c,a,b)=\int_P\mathbf1_{z+c+ax+by>0}
                    (1,\operatorname{sgn}x,\operatorname{sgn}y)^{\mathsf T}.
\]
Its derivative is exactly $A_D$ for the projected section polygon,
which contains the origin. Lemma~\ref{gen:section-matrix} makes all
principal minors positive on $\Omega$. The Gale--Nikaid\^o univalence theorem
\cite[Theorem~4]{GaleNikaido}, applied to any closed rectangle
containing two proposed preimages, proves injectivity. The six affine constraints also have independent differentials:
order them by the cut most recently added. For the first cut vary its offset;
for the second vary its offset and the slope $a$ of the quartering-fiber proof.
The resulting diagonal blocks are nonsingular by the positive section mass
and $C'(a)>0$ after median reduction. For the third cut the diagonal block
is $D\Psi=A_D$. This block triangular minor has rank $1+2+3=6$.
Invertibility of the completion block gives
$C^1$ dependence on the original pair, and uniqueness glues the
local inverses. The volume formulas give semialgebraicity. This
completes Theorem~\ref{gen:completion3}.

\paragraph{The generic metric.}
Let $\mathcal N$ consist of pairs of affine pairwise triples and compatible
metrics. Locally choose ordered representatives for the three normal rows of
$A$. Formula~\eqref{eq:gen:metric-fiber} identifies the incidence space with
$\mathcal B_3(P)\times(0,\infty)^3$ in these coordinates: the inverse reads
off the diagonal of $AQA^{\mathsf T}$. Thus $\mathcal N$ is a six-dimensional
$C^1$ semialgebraic manifold, and the fiber-dimension theorem gives
\eqref{eq:gen:pairwise-exceptional}.

Consider the projection $\pi:\mathcal N\to\operatorname{Sym}_3^+$.
Its differential has a kernel precisely when a nonzero tangent
variation of the affine triple preserves orthogonality at fixed
$Q$. Since the six affine equal-volume constraints already have
rank six, this is equivalent to failure of rank nine for their
combination with the three orthogonality equations. Domain and
target of $\pi$ have the same dimension, so this is exactly its
critical locus. A smooth
constant-rank stratification of this semialgebraic critical locus
shows that its image has dimension at most five. Adjoining its
real Zariski closure gives a proper algebraic exceptional set.
Finally the median offsets are continuous in their directions, so
$\mathcal B_3(P)$ is closed in the compact space
$(\mathbb {RP}^2)^3$. Over compact metric subsets the projection
is proper. Lemma~\ref{gen:continuation} gives the locally constant count, also after the finite ordering and orientation quotients.

For a simplex, this metric variation is the usual variation of shape.
If the columns of $E$ are its edge vectors relative to one vertex,
pullback to a reference simplex preserves relative volumes and puts
$Q=(E^{\mathsf T}E)^{-1}$ on normal covectors. Matrix inversion is
a rational isomorphism of positive definite cones. Clearing a power
of the positive determinant transfers the exceptional algebraic set
to the six Gram entries. The theorem therefore concerns arbitrary
nearby tetrahedral shapes.

\subsection{Proofs of the construction and counting theorems}
\subsubsection{The planar lift}
\begin{proof}
The map from the two child totals to their sum and $L_v$ has
determinant $-1$. Iterating it and then mapping the three block
totals to $(\sum W_i,\sum W_iB_i)$ proves affine independence
and the volume formula.

Represent the barycentric coordinates by normalized independent
exponentials. Disjoint block sums are independent gamma variables.
At each node the negative child's fraction of the parent sum is
$\operatorname{Beta}(a_v,b_v)$ and independent of that sum
\cite{BobeckaWesolowski}. Induction up the tree gives independent
ratios, independent also of the three normalized block totals $W$.
Their median thresholds give fair signs. The remaining cuts, perpendicular
to every coordinate cut, depend only on $Y=\sum W_iB_i$; hence each
lifted cell has $2^{-(n-2)}$ times its planar mass.

The projected density gives zero mass to lines, so Theorem~\ref{intro:orthogonal-existence}
in dimension two gives the required quartering.
Finally the Gram form of the construction is the pullback of
$\operatorname{diag}((\ell_v^2)_v,G_\Delta)$ by one fixed
invertible edge matrix. Its cone has dimension $n+1$; removing
scale and finite vertex relabelings leaves dimension $n$.
\end{proof}

\subsubsection{The regular simplex}
\begin{proof}
Split the vertices into at most three power-of-two blocks, and use
balanced binary trees. The centered vertex Gram matrix is a positive multiple of $I-J/(n+1)$, where $J$ is the all-ones matrix. Thus the Euclidean inner product on zero-sum coefficients agrees up to scale with that of the corresponding normals. The coefficient vectors are mutually
orthogonal Haar vectors in the regular barycentric model; each is
also perpendicular to all block-constant vectors. Their signs are
independent and fair. The remaining block-constant space has
dimension at most two. Its induced Dirichlet measure has a median
in dimension one or a perpendicular quartering in dimension two,
completing the partition. For a single block no further cut is needed.
\end{proof}

\subsubsection{Axial rigidity}
\begin{proof}
Median signs vary continuously almost everywhere by Lemma~\ref{gen:cut-flux}.

We first establish a strict correlation property. Let $Z=v\cdot Y$,
$v\ne0$, and let $h$ be the median of $Z+bX$, with $b>0$.
Conditional fairness makes the conditional $Y$-law on either side
of $X=0$ equal to $\mu$. Thus
\[
 A_+:=\Pr(Z+bX>h\mid X>0)\ge\Pr(Z>h)
 \ge A_-:=\Pr(Z+bX>h\mid X<0).
\]
If the correlation $A_+-A_-$ vanished, medianity would give
$A_+=A_-=1/2$, so $h$ would also be the median of $Z$.
Choose an interior projected point with $Z=h$. Its fiber straddles
zero. Positive density near interior points on both sides then
gives positive mass with $Z<h<Z+bX$ and with $Z+bX<h<Z$,
contradicting equality. The correlation therefore has the strict
sign of $b$.

Now consider any sequence of solutions with $L\to\infty$ and
physical orthonormal normals $n_i=(a_i,v_i)$ converging along a
subsequence. Two nonzero limiting $a_i$ would give pulled-back
median signs tending to $\pm\operatorname{sgn}X$, with correlation
$\pm1$. Exactly one limit is therefore nonzero, and its absolute
value is one, since $\sum a_i^2=1$. Orient and label it as $n_0$.
The remaining $v_i$ tend to an orthonormal frame $R_i$ of
$\mathbb R^m$. Put $b_i=-La_i$ for $i\ge1$. These numbers are
bounded, since an unbounded one would again give sign
$\pm\operatorname{sgn}X$. At any finite limiting value, the other
median forms are $R_i\cdot Y-b_iX-h_i$; the strict correlation
property forces $b_i=0$. Hence every solution approaches $b=0$
in these rescaled coordinates.

Near the limiting axis use the frame chart
\begin{equation}\label{eq:gen:stretch-chart}
 \begin{gathered}
 \epsilon=L^{-2},\quad c=\sqrt{1-\epsilon|b|^2},\quad B=\sum_i b_iR_i,\\
 n_0=(c,L^{-1}B),\qquad
 n_i=(-L^{-1}b_i,R_i-\epsilon b_iB/(1+c)).
 \end{gathered}
\end{equation}
It covers all such frames: if $V$ is the matrix of nonfirst
$Y$-rows, then $VV^{\mathsf T}=I-\epsilon bb^{\mathsf T}$;
its positive square root is
$P=I-\epsilon bb^{\mathsf T}/(1+c)$, and $R=P^{-1}V$ is
orthogonal. Pulling back and dividing the first form by $Lc$ gives
reference covectors
\[
 w_0=(1,\epsilon B/c),\qquad
 w_i=(-b_i,R_i-\epsilon b_iB/(1+c)).
\]
They are independent and smooth through $\epsilon=0$. Let
$\Phi_i(\epsilon,b,R)$ be the median-reduced correlation of
signs $0,i$.

The pulled-back sections have independent normals through $\epsilon=0$;
Lemma~\ref{gen:cut-flux} gives joint $C^1$ dependence on $(\epsilon,b,R)$.
Conditional fairness gives $\Phi(\epsilon,0,R)=0$. At the limiting
frame, differentiation in $b$ yields
\[
 \partial_{b_j}\Phi_i(0,0,R)=-2\delta_{ij}D_i(R),\qquad
 D_i(R)=\int_{R_i\cdot Y=h_\mu(R_i)}\int|X|\rho(X,Y)\,dX\,dS>0.
\]
The median-offset term disappears by conditional fairness.
Thus $\partial_b\Phi$ is invertible at every $(0,0,R)$.
The parameter-dependent implicit function theorem makes $b=0$ the unique
nearby solution. Compactness of the frame space makes the neighborhood
uniform in $R$; the preceding limit argument puts every solution there
for sufficiently large $L$.

All remaining cuts now depend only on $Y$. Every moment containing
the first sign vanishes by conditional fairness; all others are
the base moments. This proves both bijections and the local germ
correspondence. Ordering the $m$ selected correlations first, the
differential has blocks
$\left(\begin{smallmatrix}T&0\\ *&J_{\rm base}\\ *&0\end{smallmatrix}\right)$
with $T$ invertible. Eliminating its columns gives the exact rank
increment $m$. Only one orthogonal cut can have the axial normal,
so deletion creates no multiplicity in the unoriented quotient.
\end{proof}

\subsubsection{Planar examples and the proof of Theorem~\ref{gen:exact-counts}}
The directions of the three medians through vertices control how the number
of quarterings changes with the metric. Metrics are taken up to scale.

\begin{lemma}[Planar seeds]\label{gen:planar-seeds}
For any three positive integer Dirichlet weights on a triangle, the metric
disk contains open regions with exactly one and exactly three perpendicular
quarterings, all simple.
\end{lemma}
\begin{proof}
First consider any triangle with positive interior polynomial
density. In angular normal coordinates every oriented median has
a unique quartering partner $B(\alpha)\in(\alpha,\alpha+\pi)$,
with $B'>0$ and $B(B(\alpha))=\alpha+\pi$. Indeed, as the second
median chord rotates, its positive overlap with the first decreases
strictly from $1/2$ to zero: the flux is the integral of
$s-\bar s$ over an initial segment of that chord, which is strictly
negative. Interchanging the two chords gives the positive derivative
in the other angle. This argument is unchanged by positive weights.

Set $m=\alpha+B(\alpha)$, $\delta=B(\alpha)-\alpha$, and
$f(m)=-\cos\delta$; then $f(m+\pi)=-f(m)$ and
$f'=\sin\delta(B'-1)/(B'+1)$. Up to scale an inverse metric is
$\left(\begin{smallmatrix}1+x&y\\y&1-x\end{smallmatrix}\right)$,
$x^2+y^2<1$. Orthogonality is the scalar equation
\begin{equation}\label{eq:gen:partner-scalar}
 F(m;x,y)=x\cos m+y\sin m-f(m)=0.
\end{equation}
Its discriminant is the envelope
$\gamma=f(\cos m,\sin m)+f'(-\sin m,\cos m)$, and
\[
 |\gamma|^2=\cos^2\delta+
       \sin^2\delta\left(\frac{B'-1}{B'+1}\right)^2<1.
\]
Its image is compact in the open disk. It is semialgebraic in the circle
coordinates: clipped polygon masses are piecewise rational in the cut
coefficients, and the partner is characterized by its quarter-mass equation. 

Exchanging the partners changes $m$ by $\pi$, so roots modulo $\pi$
count unordered quarterings. Their number is odd by antiperiodicity
and locally constant at regular parameters by Lemma~\ref{gen:continuation}.
At a rank-one boundary metric $2aa^{\mathsf T}$ the orthogonality
equation says that one partner is perpendicular to $a$; uniqueness
of the partner gives one simple root. Thus a boundary annulus has
count one.

Except for finitely many image values, a discriminant crossing is
a single fold. To check this without a genericity assumption on
the density, split the semialgebraic envelope into analytic arcs.
There $\gamma'=(f+f'')(-\sin m,\cos m)$. Zero intervals of
$f+f''$ have constant images. Two distinct regular parameters
modulo $\pi$ have linearly independent tangent directions, so
their coincidences are isolated and, by semialgebraicity, finite.
At a regular image outside these exceptional values there is one double
root: $F_{mm}=-(f+f'')\ne0$, and variation of the metric in the normal
direction unfolds it. Thus it suffices to prove that $\gamma$ is
nonconstant. Indeed, choose a generic path from the one-root annulus
to a regular arc and across it, avoiding the finite exceptional values.
At its first transverse crossing the count changes by two. Since every
metric has a quartering, the new count is three.

Suppose instead that $\gamma$ were constant. Its value would be an
interior metric for which $F(m)=0$ for every $m$: every median would
quarter the mass with its perpendicular median. Fix this metric.
In unoriented normal directions perpendicularity is a fixed-point-free
involution. There are exactly three medians passing through a vertex,
one per vertex: sweeping a ray from one incident edge to the other
increases the swept mass strictly from zero to one. These three lines
are distinct since an edge cannot bisect the triangle. Consequently
one vertex median has a perpendicular median containing no vertex.

Write their normals as $u(\theta)$ and $u(\theta+\pi/2)$, and put the
vertex at the origin. Choose orthonormal coordinates for the fixed metric and take unit normals.
The second median has analytic offset near this pair: its clipped polygon has fixed incidences, the density is polynomial,
and its section mass is positive. Take its half-plane omitting the
origin. Within that half-plane the first cut bisects the mass, since the
pair quarters. Its offset $d(\theta)$ is therefore analytic as well:
the clipped polygon has fixed incidences there, the intersection of the
two medians is interior, and the first cut has positive section mass
on this side. Yet an analytic offset is impossible at the vertex.

For the Dirichlet density, the obstruction is explicit. Use barycentric
edge coordinates $x,y\ge0$ at the vertex. Locally the density is $x^{i-1}y^{j-1}g(x,y)$ with
$i,j\ge1$ integers and $g$ polynomial, $g(0,0)>0$. The median line reads
$a x-b y=d$, with $a,b>0$ analytic in $\theta$. At the vertex direction
$d=0$, and $d'\ne0$: by \eqref{eq:gen:median-derivative}, $d'$ is the section mean of
$u'\cdot X$, which has one strict sign along the open median chord
issuing from the origin. Away from a small corner the cut mass is analytic.
In the corner, integrating first in $y$ up to $(ax-d)/b$ gives an analytic
continuation from $d<0$. For $d>0$ the interval $0<x<d/a$ must be removed
from that continuation. Its signed contribution, by $x=dt/a$, is
\[
 \frac{(-1)^j g(0,0)}{j a^i b^j}
 \left(\int_0^1 t^{i-1}(1-t)^j\,dt\right)d^{i+j}
 +O(d^{i+j+1}).
\]
Here a fixed small corner can be cut off by $x+y=\varepsilon$; its other
intersection with the median is transverse and varies analytically.
The coefficient displayed above is nonzero. Since $d$ changes sign
transversely, the analytic continuation of the constant half-mass from
one side cannot equal the half-mass on the other. This contradiction
proves the lemma.
\end{proof}

\begin{proof}[Proof of Theorem~\ref{gen:exact-counts}]
Corollary~\ref{gen:integer-split} reduces both counts to seeds of total
weight $n+1$. For one partition take the interval of weights $(1,n)$:
its median is unique and the initial frame dimension is zero. For three
partitions in dimension $n\ge2$, take a triangle with weights
$(1,1,n-1)$. Its total weight is $n+1$, and Lemma~\ref{gen:planar-seeds}
supplies a metric with exactly three simple quarterings.

Split to unit weights in each case. Every partition descends, so all
remain full. A split in dimension $m$ adds $m$ to the reduced rank;
starting with rank zero on the interval or one on the triangle gives
$n(n-1)/2$. In dimension three the split is symmetric and its mirror
is the prescribed cut of Theorem~\ref{thm:c11-edge}. Sufficient stretch
is essential: splitting the projected midpoint of an edge back to its
original endpoints in a regular tetrahedron also gives partitions
omitting this mirror, by Theorem~\ref{intro:circles}.
\end{proof}

\Needspace{7\baselineskip}
\section{Orthogonal partitions of tetrahedra}
\label{sec:tetrahedra}

The centered partitions are determined by the three bimedians, and their
fullness by one further sign moment. We first give this classification,
then prove that its continuous families are isolated. The final results
reduce partitions with a face-parallel or edge-containing cut to planar
quadrisections and exclude three strict singleton cuts in the regular case.

\subsection{Definitions and main results}
\label{subsec:tetra-results}

Let $T\subset\mathbb R^3$ be a nondegenerate tetrahedron, with the
configuration and counting conventions of Section~\ref{gen:definitions}.
For a pairwise partition, \eqref{eq:gen:triple-moment} expresses all eight
volume fractions through the single moment $\tau$; fullness is $\tau=0$.

A partition is \emph{centered} if all three planes pass through
the centroid $g$. A \emph{bimedian} joins the midpoints of opposite
edges. Use the reference tetrahedron
\[
 S_0=\operatorname{conv}\{(1,1,1),(1,-1,-1),(-1,1,-1),(-1,-1,1)\}
\]
and write $T=g+AS_0$. Set $b_i=Ae_i$ and $B=A^{\mathsf T}A$;
the vectors $2b_i$ are the three bimedians. The first two results
classify the complete centered solution space.

\begin{theorem}\label{c11:centered-pairwise}
A centered orthogonal partition of $T$ exists if and only if one
$b_i$ is perpendicular to the other two. Equivalently, in at least two of the three pairs of opposite edges, the two edges have equal length. The complete centered
solution set is described by the following table.
\begin{center}
\begin{tabular}{@{}ll@{}}
\hline
Condition on $B$ & Centered partitions\\
\hline
No coordinate axis is an eigenline & None\\
At least one coordinate eigenline, except below & One\\
$B$ diagonal with exactly two equal entries & One circle\\
$B$ scalar & A bouquet of three circles\\
\hline
\end{tabular}
\end{center}
In particular, a nonconstant centered family exists precisely when
the three bimedians are mutually perpendicular and at least two
have equal lengths. Up to isometry, these tetrahedra are
\[
 \operatorname{conv}\{(a,0,0),(-a,0,0),(0,a,h),(0,-a,h)\},
 \qquad a,h>0.
\]
\end{theorem}

\begin{theorem}\label{c11:centered-full}
A tetrahedron has a centered full orthogonal octasection if and
only if four edges forming a four-cycle have equal lengths.
Equivalently, up to isometry and relabeling, it has the form
\begin{equation}\label{c11:full-family}
 T(a,b,h)=\operatorname{conv}
 \{(a,0,0),(-a,0,0),(0,b,h),(0,-b,h)\},
 \qquad a,b,h>0.
\end{equation}
The centered full-octasection count is one for every nonregular
tetrahedron in this family, three for a regular tetrahedron, and
zero for every other tetrahedron.
\end{theorem}

\begin{theorem}[Isolation of centered families]
\label{thm:v33:centered-positive-components}
If the centered solution set of a tetrahedron is infinite, it is an
isolated connected component of the space of all orthogonal partitions.
\end{theorem}

By Theorem~\ref{c11:centered-pairwise}, this component is a circle for a
nonregular tetrahedron and three circles meeting at one point for a regular one.

\begin{figure}[htbp]
 \centering
 \includegraphics[width=.68\linewidth]{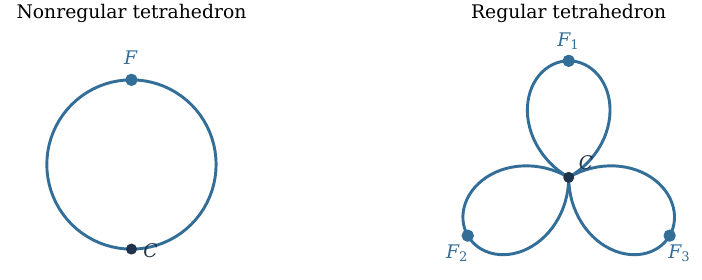}
 \caption{The centered solution space when it contains a nonconstant
 family (topological schematics). The marked points $F$ and $F_i$
 represent full octasections. At the regular shape, three circles
 meet at the coordinate partition $C$. By
 Theorem~\ref{thm:v33:centered-positive-components}, each displayed
 space is an isolated connected component of the entire pairwise
 configuration space of the fixed tetrahedron.}
 \label{fig:centered-components}
\end{figure}

Figure~\ref{fig:centered-components} shows the continuous centered families.
Each point represents a partition, with the full ones marked $F$ or $F_i$.
The common point $C$ represents the coordinate planes, whose eight cell
volumes are unequal.

A face-parallel cut is parallel to one face; an edge-containing cut
contains an entire edge. In the following face-parallel statement,
$V$ is a vertex, $B$ is its opposite face, and $p$ is the orthogonal
projection of $V$ onto the plane of $B$.

\begin{theorem}\label{thm:c11-face}
Orthogonal partitions with a cutting plane parallel to $B$ are in
bijection with perpendicular equal-area quadrisections of $B$ whose
two lines meet at $p$. Every such partition is a full octasection.
\end{theorem}

\begin{theorem}\label{thm:c11-edge}
Suppose one plane of an orthogonal partition of $ABCD$ contains the
edge $AB$. It is then the Euclidean reflection plane fixing $A,B$ and
exchanging $C,D$, and the partition is a full octasection.
Conversely, every reflection plane of a tetrahedron belongs to a full
orthogonal octasection. Those containing a prescribed reflection plane
are in bijection with perpendicular quadrisections of a triangle with
density proportional to one barycentric coordinate.
\end{theorem}

A \emph{strict singleton cut} separates one vertex from the other three
and passes through no vertex.

\begin{theorem}[Exclusion of three strict singleton cuts]
\label{thm:v24:regular-s3-exclusion}
For a regular tetrahedron, suppose three mutually orthogonal
median planes are strict singleton cuts and every pair quarters
volume. Orient each cut positively at its isolated vertex.
These vertices are distinct, and the common positive cell has
volume fraction strictly less than $1/8$. Thus no such triple
is a full octasection.
\end{theorem}

Theorem~\ref{c11:centered-full} gives Theorem~\ref{intro:center}; the
classification of noncentral full octasections remains open.
\FloatBarrier

\subsection{Centered medians and centered partitions}
\label{subsec:tetra-centered}
We use two symmetries of $S_0$: even coordinate sign changes and coordinate
permutations preserve it. Averaging its volume measure with its central
reflection makes all eight coordinate signs equiprobable at fixed absolute
coordinates $|X|=(|X_1|,|X_2|,|X_3|)$.

\begin{lemma}[Centered medians]\label{c11:median}
A plane through the centroid of a tetrahedron bisects its volume if
and only if it contains the line joining the midpoints of some pair
of opposite edges. In particular, a centered plane with normal
$q=(q_1,q_2,q_3)$ bisects $S_0$ if and only if $q_1q_2q_3=0$.
\end{lemma}
\begin{proof}
If one coordinate of $q$ vanishes, changing the signs of the other two
coordinates reverses $q\cdot X$ and preserves $S_0$, so the plane is median.
For the converse, reverse the normal and change two coordinate signs if
necessary, so all $q_i>0$.

The central octahedron $|x_1|+|x_2|+|x_3|\le1$ contributes zero to the
integral of $\operatorname{sgn}(q\cdot X)$. The rest of $S_0$ consists
of four reflected copies of
\[
 D=\{a_i>0:\ 1<a_1+a_2+a_3<1+2\min_i a_i\},
\]
with coordinate sign product $+1$. For positive $A_i$, away from the
triangle-equality cases,
\[
 \sum_{\varepsilon_1\varepsilon_2\varepsilon_3=1}
 \operatorname{sgn}(\varepsilon_1A_1+\varepsilon_2A_2+\varepsilon_3A_3)
 =-2\mathbf1_{\,2\max A_i<A_1+A_2+A_3}.
\]
Apply this with $A_i=q_i a_i$. The sum is negative on an open subset
of $D$: a positive rescaling of $(q_1^{-1},q_2^{-1},q_3^{-1})$ lies in
$D$ and makes the three $A_i$ equal. Hence the sign integral is strictly
negative. Thus the plane is median exactly when some $q_i=0$.
The three coordinate axes join opposite-edge midpoints of $S_0$;
affine pullback gives the stated criterion for $T$.
\end{proof}

\begin{lemma}[Pairwise reduction]\label{c11:correlation}
Two centered median planes of $S_0$, with nonzero normals $q,r$,
quarter its volume if and only if $q\cdot r=0$.
\end{lemma}
\begin{proof}
By Lemma~\ref{c11:median}, each normal has a zero coordinate. If the
zeros can be chosen in the same position, project onto the complementary
coordinate plane. The projected law is invariant under a quarter-turn:
for the $xy$-projection, $(x,y,z)\mapsto(-y,x,-z)$ permutes the vertices
of $S_0$. Each quadrant therefore has mass $1/4$, and a sector from
any fixed ray has strictly increasing mass with its angle. The quartering
partner is precisely the perpendicular line.

Otherwise we may write $q=(a,b,0)$ and $r=(c,0,d)$ after a coordinate
permutation. Their sign product is unchanged by $X\mapsto-X$, so the symmetrized measure
gives the same integral. Averaging its eight sign choices at fixed $|X|$ gives
\[
 \mathbb E[\operatorname{sgn}(q\cdot X)\operatorname{sgn}(r\cdot X)\mid |X|]
 =\operatorname{sgn}(ac)\,
 \mathbf1_{|aX_1|>|bX_2|}\mathbf1_{|cX_1|>|dX_3|}.
\]
For $ac\ne0$, the indicated set has positive measure because the
support contains a neighborhood of the origin. The expectation
vanishes precisely when $ac=q\cdot r=0$. Since both planes are
medians, zero sign correlation is equivalent to four equal volumes.
\end{proof}

\begin{proof}[Proof of Theorem~\ref{c11:centered-pairwise}]
Pull the physical normals back by $q_i=A^{\mathsf T}u_i$.
The preceding lemmas say that the normalized $q_i$ form an orthonormal
basis, each vector having a zero coordinate. Physical orthogonality is
$q_i^{\mathsf T}B^{-1}q_j=0$ for $i\ne j$, so this basis must diagonalize
$B$. These conditions are also sufficient.

Such a basis contains a coordinate axis. Otherwise each column of its
orthogonal matrix has exactly one zero. Two zeros in one row force the
third column onto that axis; zeros in different rows leave two columns
with only one common nonzero coordinate, contradicting orthogonality.
Once an axis is chosen, the remaining two vectors diagonalize the
complementary $2\times2$ block. Its eigenlines are unique unless the block
is scalar, when every rotation is allowed. A coincidence between the
axial eigenvalue and a block eigenvalue does not change this argument.
This gives the four cases of the table. When $B$ is scalar, the three
circles meet only at the coordinate frame.

For example, the circle with fixed axis $e_3$ in $S_0$ is
\begin{equation}\label{c11:circle}
 q_1=(\cos\theta,\sin\theta,0),\quad
 q_2=(-\sin\theta,\cos\theta,0),\quad q_3=e_3,
 \qquad \theta\in\mathbb R/(\tfrac\pi2\mathbb Z).
\end{equation}
The period follows from the unordered, unoriented convention.

For the vertex ordering used in $S_0$, the differences of two pairs
of opposite squared edge lengths are $16b_1\cdot b_2$ and
$16b_1\cdot b_3$. Relabeling proves the edge-length criterion.
Finally, diagonal $B$ with $|b_2|=|b_3|$ gives perpendicular,
equal-length vectors $b_2+b_3$ and $b_2-b_3$, both perpendicular
to $b_1$. These are the half-edge vectors of the displayed family,
with $a=|b_2+b_3|$ and $h=2|b_1|$.
\end{proof}

\begin{proof}[Proof of Theorem~\ref{c11:centered-full}]
For the circle~\eqref{c11:circle}, integrate $\operatorname{sgn}z$ along
the vertical fiber, whose endpoints are $|x+y|-1$ and $1-|x-y|$.
Write the resulting signed length as $W(x,y)$. On $0<x,y<1$ the fiber
has positive length and positive midpoint $\min(x,y)$, hence $W>0$.
It is symmetric in $x,y$ and changes sign under a quarter-turn.
The product of the first two cut signs changes sign under the same turn,
so all four quadrants contribute equally to the triple moment $\tau$.

In the first quadrant that product is $\operatorname{sgn}(\phi-\theta)$,
where $\phi=\arg(x+iy)$ and $0\le\theta\le\pi/2$. Thus $\tau$ is a positive
multiple of the difference between the $W$-masses above and below the ray
$\phi=\theta$. It decreases strictly with $\theta$ and vanishes at
$\pi/4$ by symmetry. Each circle therefore contains exactly one full partition.
The three normal systems are
\begin{equation}\label{c11:three-full}
 e_i,\quad (e_j+e_k)/\sqrt2,\quad(e_j-e_k)/\sqrt2,
 \qquad \{i,j,k\}=\{1,2,3\}.
\end{equation}

Pull a centered full octasection of $T$ back to $S_0$. By
Lemma~\ref{c11:correlation}, its pulled-back normals are perpendicular,
so it is one of~\eqref{c11:three-full}. Its physical normals are
perpendicular precisely when this basis diagonalizes $B$, namely
\begin{equation}\label{c11:full-metric}
 b_i\cdot b_j=b_i\cdot b_k=0,\qquad |b_j|=|b_k|.
\end{equation}
The opposite edges with midpoints $g\pm b_i$ have half-edge vectors
$b_j+b_k$ and $b_j-b_k$. Conditions~\eqref{c11:full-metric} make
these vectors perpendicular to each other and to $b_i$, yielding
\eqref{c11:full-family}. For each $i$ there is only one such normal
system. If two indices satisfy~\eqref{c11:full-metric}, all three
$b_i$ are perpendicular and have equal length; hence $B$ is scalar
and $T$ is regular. This proves the count.

Finally, write opposite edges as $M\pm u$ and $N\pm v$, with
$w=N-M$. The four cross-edge squared lengths $|w\pm u\pm v|^2$
are equal exactly when
$w\cdot u=w\cdot v=u\cdot v=0$. This is precisely the coordinate
form~\eqref{c11:full-family}.
\end{proof}

For example, $T(1,1,1)$ has edge lengths
$2,2,\sqrt3,\sqrt3,\sqrt3,\sqrt3$ and one centered full octasection.
The family~\eqref{c11:full-family} is regular exactly when
$a=b$ and $h=\sqrt2a$.

\subsection{Isolation of centered components}
\label{subsec:tetra-centered-components}

The isolation proof has three local cases: an invertible transverse differential,
the common point of the regular bouquet, and one exceptional axial frame.
We treat them in this order; compactness then gives the global statement.

By Theorem~\ref{c11:centered-pairwise}, a tetrahedron with a nonconstant
centered family is similar to
\[
 T_\lambda=\operatorname{diag}(1,1,\lambda)S_0,\qquad \lambda>0.
\]
All coordinates below are physical, so
$\operatorname{vol}(T_\lambda)=8\lambda/3$.
For a physical orthonormal frame $U=(u_1,u_2,u_3)$, put
$m_i=h_{T_\lambda}(u_i)$ and
\[
 \mathbf F_\lambda(U)=(F_{12},F_{13},F_{23}),\qquad
 F_{ij}=\mathbb E_{T_\lambda}
   [\operatorname{sgn}(u_i\cdot X-m_i)
    \operatorname{sgn}(u_j\cdot X-m_j)].
\]
These are the median-reduced moments of~\eqref{eq:gen:pair-moment}.
The centered circle $\mathcal C_\lambda$ is represented by
$U_\theta=R_z(\theta)$, with $\theta$ modulo $\pi/2$. Reflection symmetries reduce the differential calculation to
$0\le\theta\le\pi/4$; their images cover the whole circle.
Pulling the cuts back to $S_0$ and applying Lemma~\ref{gen:cut-flux} gives joint
$C^1$ dependence on frame and stretch, also at incidences.

\begin{lemma}[Transverse differential]
\label{thm:v33:axial-circle-differential}
Use $U=U_\theta\exp[w]_\times$, where $[w]_\times v=w\times v$, and put
$t=\tan\theta\in[0,1]$. Then
\begin{equation}\label{eq:v33:axial-circle-jacobian}
 \begin{gathered}
 D\mathbf F_\lambda(U_\theta)=
 \begin{pmatrix}0&0&0\\-A_\lambda&-B_\lambda&0\\B_\lambda&-A_\lambda&0\end{pmatrix},
 \qquad
 A_\lambda(t)=\frac{t(1-t)}{2\lambda(1+t)^2\sqrt{1+t^2}},\\
 B_\lambda(0)=\frac{\lambda-\lambda^{-1}}2,
 \qquad B_\lambda(1)=\frac{4\lambda^2-3}{4\lambda\sqrt2}.
 \end{gathered}
\end{equation}
The transverse block is invertible except at
$(t,\lambda)=(0,1)$ and $(1,\sqrt3/2)$, where the whole differential vanishes.
\end{lemma}
\begin{proof}
Symmetry determines the matrix form. The axial half-turn makes $F_{12}$
even in the transverse variables, and the third column is tangent to the
centered circle. For $R=R_z(\pi/2)$ and
$J=\left(\begin{smallmatrix}0&-1\\1&0\end{smallmatrix}\right)$,
$RT_\lambda=-T_\lambda$ and reversal of all signs preserves pair moments.
With $G=(F_{13},F_{23})$, relabeling the first two cuts gives
\[
 G(U_\theta e^{[Rw]_\times})=J G(U_\theta e^{[w]_\times}).
\]
Thus the transverse block commutes with $J$ and is determined by its first row.

Write $L_1=x+ty$, $L_2=-tx+y$, $L_3=z$. The half-section
$D_+=\{|x|+|y|\le1,\ x+ty>0\}$ has (unnormalized) moments
\[
 I_x=\frac{2t+1}{3(1+t)^2},\qquad
 I_y=\frac{t(t+2)}{3(1+t)^2}.
\]
These follow by triangulating from the origin its vertices
$(1,0),(0,1),(-t/(1+t),1/(1+t)),(t/(1+t),-1/(1+t))$.
A $w_1$ rotation fixes the first cut and changes the third form by
$(tx-y)w_1/\sqrt{1+t^2}$, whose median derivative is zero. Lemma~\ref{gen:cut-flux} gives
$\partial_{w_1}F_{13}=4(tI_x-I_y)/(\vol(T_\lambda)\sqrt{1+t^2})=-A_\lambda$.

Since $A_\lambda>0$ for $0<t<1$, only the endpoint values of $B_\lambda$
are needed. Under $w_2=\varepsilon\sqrt{1+t^2}$ the forms change to
$L_1-\varepsilon(1+t^2)z$ and $z+\varepsilon L_1$.
The section $L_1=0$ has coarea measure $dy\,dz$ and bounds
$|y|\le1$, $\lambda[-1+(1-t)|y|]\le z\le\lambda[1-(1+t)|y|]$.
Let $\bar z$ be its mean, $H=I_x+tI_y$, and
$J_t=\int_{L_1=0,z>0}(1+t^2)(z-\bar z)\,dy\,dz$.
Direct integration at the endpoints gives
\[
\begin{array}{c|ccc}
 t&\bar z&H&J_t\\ \hline
 0&0&1/3&\lambda^2/3\\
 1&-\lambda/3&1/2&2\lambda^2/3
\end{array}
\]
The two moving sections contribute
$\partial_{w_2}F_{13}=4(H-J_t)/(\vol(T_\lambda)\sqrt{1+t^2})$.
This gives the stated endpoint values. The determinant is
$A_\lambda^2+B_\lambda^2$, proving the rank assertion.
\end{proof}

\begin{lemma}[The regular bouquet at its common point]
\label{lem:v33:regular-coordinate-germ}
At the coordinate partition of $S_0$, in rotation-vector coordinates,
\begin{equation}\label{eq:v33:regular-coordinate-germ}
 \mathbf F_1(e^{[w]_\times})
 =-\tfrac12(w_1w_2,w_1w_3,w_2w_3)+O(\|w\|^3),
\end{equation}
with derivative remainder $O(\|w\|^2)$. Its local zero set is exactly
the three coordinate rotation axes.
\end{lemma}

\begin{proof}
Each median is analytic near a coordinate normal: all vertex heights
stay away from zero, the median cuts have a fixed two-positive,
two-negative type, and their cap density is positive.
For the pair $x\ge0,y\ge0$, the only nonsimple clipped vertices are
$(0,0,\pm1)$. At each, two cuts and two tetrahedral faces meet,
and every three of the four normals are independent. Choose three
planes as affine coordinates $z_1,z_2,z_3$. The fourth has equation
$a_1z_1+a_2z_2+a_3z_3+\eta=0$, with $a_1a_2a_3\ne0$.
The tetrahedron exchanged between its two analytic triangulations
has vertices $0,-(\eta/a_1)e_1,-(\eta/a_2)e_2,-(\eta/a_3)e_3$;
its volume is an analytic coefficient times $\eta^3$.
Consequently the piecewise analytic volume is $C^2$, with a common
quadratic jet and uniform cubic value and quadratic derivative
remainders. The two vertex neighborhoods are disjoint, so simultaneous
changes give the same conclusion. The other pairs are identical by
coordinate permutation.

The map vanishes on all coordinate rotation axes. The coordinate
half-turn characters therefore force its quadratic jet to be
$(c_{12}w_1w_2,c_{13}w_1w_3,c_{23}w_2w_3)$.
At $\lambda=1$, the preceding differential gives
$\partial_{w_1}F_{13}|_{R_z(\theta)}=-\theta/2+O(\theta^2)$.
Right-rotation and rotation-vector coordinates differ by $I+O(\theta)$;
the differential is already $O(\theta)$, so this gives $c_{13}=-1/2$.
Cyclic coordinate permutations give the other two coefficients.

Choose $|w_3|=\max_i|w_i|$, up to permutation, and put $v=(w_1,w_2)$.
Along $(tv,w_3)$ the transverse derivative of $(F_{13},F_{23})$ is
$-w_3 I/2+O(w_3^2)$ uniformly. The map vanishes at $(0,0,w_3)$;
integrating this derivative along the segment gives
$\|(F_{13},F_{23})\|\ge |w_3|\|v\|/4$ for small $w$.
Thus every zero lies on an axis, and every axis consists of zeros.
\end{proof}

\subsection{The critical axial configuration}
\label{sec:v33:critical-germ}

The only nonregular exceptional value in the transverse differential
is $\lambda=\sqrt3/2$, at the reflection frame $\theta=\pi/4$.
Although its differential vanishes, no further solution branch
passes through it.

\begin{lemma}[Rigidity of the critical axial frame]
\label{thm:critical-axial-rigidity}
For $T_{\sqrt3/2}$, the local orthogonal pairwise partition set at
$U_{\pi/4}$ is exactly the centered circle $\mathcal C_{\sqrt3/2}$.
Already the two equations $F_{13}=F_{23}=0$ define this circle locally.
\end{lemma}

\begin{proof}
Rotate the horizontal coordinates and pull back by
$\operatorname{diag}(\sqrt2,\sqrt2,\lambda)$ to the reference tetrahedron
\begin{equation}\label{eq:v33:critical-K}
 \begin{aligned}
 K&=\operatorname{conv}\{(\pm1,0,1),(0,\pm1,-1)\}\\
  &=\left\{ |X|\le\frac{1+Z}{2},\ |Y|\le\frac{1-Z}{2},\ |Z|\le1\right\}.
 \end{aligned}
\end{equation}
Its normalized volume density is $3/4$.  The reference covectors
are orthogonal for $\operatorname{diag}(1,1,\gamma)$, where
$\gamma=2/\lambda^2=8/3$.  Use the local chart
\[
 \ell_1=X+pY+qZ,\quad
 \ell_2=rX+Y+sZ,\quad
 \ell_3=uX+vY+Z,
\]
with
\begin{equation}\label{eq:v33:critical-chart}
 r=-p-\gamma qs,\qquad
 u=\frac{\gamma(ps-q)}{1-pr},\qquad
 v=\frac{\gamma(rq-s)}{1-pr}.
\end{equation}
Normalization by the leading coefficient is valid nearby, and solving
the three orthogonality equations gives this chart. Its centered circle is $q=s=0$.
Write $\rho=|p|+\sqrt{q^2+s^2}$.
Throughout this proof, $R=O_{C^1}(\rho^3)$ means
$|R|\le C\rho^3$ and $\|DR\|\le C\rho^2$, uniformly across
all incidence sectors.

Write the second-order median correction as
\begin{equation}\label{eq:v33:critical-mu}
 \mathcal M(p,q)=\mathbb E[(2q/3-p\zeta)|2q/3-p\zeta|],\qquad \zeta\sim\mathrm{Unif}[-1,1].
\end{equation}
It is $C^{1,1}$ and even in $p$. We claim that the
three medians satisfy
\begin{equation}\label{eq:v33:critical-medians}
 \begin{aligned}
 m_1&=-q/3+\mathcal M(p,q)+O_{C^1}(\rho^3),\\
 m_2&=s/3-\mathcal M(p,s)+O_{C^1}(\rho^3),\\
 m_3&=O(\rho^2),\qquad Dm_3=O(\rho).
 \end{aligned}
\end{equation}
At $p=q=s=0$, each of the first two planes contains an entire edge of $K$,
so an ordinary analytic Taylor expansion would be unjustified.
To prove the first formula, put $Z=2t-1$, $Y=(1-t)y$, and
$a=py-q-m$, $b=2q-py$.  For fixed $y$ the unclipped sign integral over $-t\le X\le t$ is
$2(a+bt)$; clipping changes it only near $t=0$. With the Jacobian
weight $1-t$, its positive saturation correction, after removing the
factor two and writing $a_+=\max\{a,0\}$, is
\[
 J(a,b)=\int_0^{a_+/(1-b)}(1-t)(a_+-(1-b)t)\,dt
 =\frac{a_+^2}{2(1-b)}-\frac{a_+^3}{6(1-b)^2}.
\]
For small parameters the upper limit is below one. Including negative
saturation gives $J(a,b)-J(-a,-b)=a|a|/2+O_{C^1}((|a|+|b|)^3)$.
Integrating the unclipped term gives $-q-3m$; the exact sign mean is
this quantity minus $3\int_{-1}^1[J(a,b)-J(-a,-b)]\,dy$.
The median therefore has linear term $-q/3$. Substitution of that term
in $a|a|$ gives $\mathcal M(p,q)$ with cubic error and quadratic derivative
error, by the Lipschitz derivative of $a|a|$ and implicit differentiation.
The symmetry $(X,Y,Z)\mapsto(Y,X,-Z)$ gives $m_2$: substituting
$r=-p+O(\rho^2)$ in the even $C^{1,1}$ function $\mathcal M$ preserves the remainder.
For the analytic third median, \eqref{eq:gen:median-derivative} and the
zero centroid of the square section $Z=0$ give the estimates for $m_3$.

To compute the transverse moments, put
\[
 Q(p,q,m;u,v,n)=\mathbb E_K\bigl[
 \operatorname{sgn}(X+pY+qZ-m)\,
 \mathbf1_{Z+uX+vY>n}\bigr].
\]
This function is analytic near zero. Indeed, both the upper half-tetrahedron
and its intersection with $X\ge0$ have simple vertices there, so their
clipped volumes are analytic in the plane coefficients.
At $u=v=n=0$, horizontal integration and then Lemma~\ref{gen:cut-flux} on
$[-1/2,1/2]^2=K\cap\{Z=0\}$ give
\[
 Q=\frac q4-\frac{3m}4,\qquad
 Q_u=\frac3{16}+O((|p|+|m|)^2),\qquad
 Q_v=\frac p8,\qquad Q_n=\frac{3m}2.
\]
For example, $Q_v=(3/4)\int_{-1/2}^{1/2}2Y(pY-m)\,dY=p/8$.
At $p=q=m=n=0$, reflection in $X$ makes $Q$ odd in $u$, and reflection
in $Y$ makes it even in $v$. Thus its quadratic terms in $u,v$ vanish.
Since $n=m_3=O(\rho^2)$ and $F_{13}=2Q$ when $m=m_1$, the expansion reduces to
\[
 F_{13}=\frac q2-\frac{3m_1}2+\frac{3u}8+\frac{pv}4
          +O_{C^1}(\rho^3).
\]
Substituting~\eqref{eq:v33:critical-medians} and expanding~\eqref{eq:v33:critical-chart} as
$u=-\gamma q+\gamma ps+O(\rho^3)$,
$v=-\gamma s-\gamma pq+O(\rho^3)$ and $\gamma=8/3$ now gives
\begin{equation}\label{eq:v33:critical-transverse}
 \begin{aligned}
 (F_{13},F_{23})&=G+O_{C^1}(\rho^3),\\
 G&=\bigl(-\tfrac32\mathcal M(p,q)+ps/3,
          -\tfrac32\mathcal M(p,s)-pq/3\bigr).
 \end{aligned}
\end{equation}
The second component follows by the symmetry $(X,Y,Z)\mapsto(Y,X,-Z)$.

The coercivity of the transverse term follows directly from the integral.
For $q\ge0$ its derivative is
$\partial_q\mathcal M=(4/3)\mathbb E|2q/3-p\zeta|$.
By Jensen the expectation is at least $2q/3$; by symmetry and convexity it
is minimized, as a function of the shift, at zero, where it equals $|p|/2$.
Integrating from $0$ to $q$, and using oddness for $q<0$, gives
\begin{equation}\label{eq:v33:critical-coercivity}
 q\mathcal M(p,q)\ge q^2\max\{2|p|/3,4|q|/9\}.
\end{equation}
Put $\xi=(q,s)$. The mixed terms in $\xi\cdot G$ cancel.
Since $|q|^3+|s|^3\ge\|\xi\|^3/\sqrt2$, \eqref{eq:v33:critical-coercivity} gives
$\xi\cdot G\le-c\rho\|\xi\|^2$ for an absolute $c>0$.
Both $F=(F_{13},F_{23})$ and $G$ vanish at $\xi=0$; integrate
$D_\xi(F-G)=O(\rho^2)$ along $(p,t\xi)$ to get
$\xi\cdot F\le(-c\rho+O(\rho^2))\|\xi\|^2$.
Every sufficiently small zero therefore has $\xi=0$, exactly the centered circle.
\end{proof}

\begin{proof}[Proof of Theorem~\ref{thm:v33:centered-positive-components}]
Symmetry reduces the circle parameter to $0\le\theta\le\pi/4$.
Where the transverse determinant in~\eqref{eq:v33:axial-circle-jacobian} is positive,
the equations
$(F_{13},F_{23})=0$ determine a unique local graph by the implicit
function theorem. The known circle supplies that graph and satisfies
$F_{12}=0$. There are just two exceptions: Lemma~\ref{lem:v33:regular-coordinate-germ}
gives the three axes at the regular common point, and
Lemma~\ref{thm:critical-axial-rigidity} gives the single circle at the
critical axial frame. The resulting centered locus is relatively open
by these local descriptions, closed by compactness, and connected.
It is therefore one whole isolated component.
\end{proof}

\subsection{Face-parallel and edge-containing cuts}
\label{subsec:tetra-incidences}

\begin{proof}[Proof of Theorem~\ref{thm:c11-face}]
The parallel median cuts off a cap homothetic to the whole tetrahedron
about $V$, with ratio $\lambda=2^{-1/3}$. If a second median has normal
$u$ parallel to $B$, its median offset on the cap is
$(1-\lambda)\langle u,V\rangle+\lambda h_T(u)$.
Pairwise quartering requires this to equal $h_T(u)$, so
$h_T(u)=\langle u,V\rangle$. Both remaining planes therefore contain
the altitude $Vp$. Their traces on $B$ are perpendicular lines through
$p$, and the four cones above the resulting face regions have equal
volume exactly when those regions have equal area. The parallel median
halves each cone. Conversely, such a face quadrisection gives four equal cones, each halved
by the parallel median, and hence a full partition.
\end{proof}

In particular, every nondegenerate triangle occurs as the base of a
one-parameter family of tetrahedra with a full octasection: choose any
of its Bernoulli quadrisections and place the apex at arbitrary height
above the intersection of its lines. With base height zero and apex
height $H$, the parallel median is
\[
 z=H(1-2^{-1/3}).
\]
Its common point with the other cuts differs from the centroid, whose
height is $H/4$.

\begin{proof}[Proof of Theorem~\ref{thm:c11-edge}]
A median plane through $AB$ meets $CD$ at its midpoint; hence its
barycentric equation is $\lambda_C=\lambda_D$. Let $S$ be the
volume-preserving affine involution exchanging $C,D$ and fixing $A,B$.
Write another cutting plane as $L=0$, where
$L=\sum k_i\lambda_i$. If $k_C>k_D$, then
\[
 L(X)-L(SX)=(k_C-k_D)(\lambda_C-\lambda_D)>0
 \quad\text{on }\lambda_C>\lambda_D.
\]
On this half-tetrahedron both $L$ and $L\circ S$ bisect volume:
the first by pairwise balance, the second because $S$ exchanges the halves.
Lemma~\ref{gen:ordering} excludes $k_C>k_D$, and exchanging $C,D$ excludes
$k_D>k_C$. Thus every other cut is parallel to $CD$. Orthogonality now makes the first normal parallel
to $CD$, so the first plane is its perpendicular bisector and is a
reflection plane of the tetrahedron. It halves each of the four regions
formed by the other two planes.

For the converse, write the reflected pair as $M\pm bn$, with $M,A,B$
in the reflection plane $P$. Projection onto $P$ has image
$\Delta=\conv\{M,A,B\}$. Its fiber above $Y$ is
$\{Y+tn:|t|\le b\lambda_M(Y)\}$. The projected probability measure is
\begin{equation}\label{eq:c11-weight}
 d\mu(Y)=\frac{3\lambda_M(Y)}{\operatorname{area}(\Delta)}\,dY.
\end{equation}
This is the $(2,1,1)$ case of~\eqref{eq:gen:dirichlet-density}. Apply
Theorem~\ref{intro:orthogonal-existence} in the mirror plane to this density.
Each quartering lifts to four quarter-volume cells, each bisected by reflection.
Conversely every partition containing $P$ projects to such a quartering.

\end{proof}

If two cutting planes contain edges, those edges are opposite and the
partition is centered. Indeed, after the first cut contains $AB$, any
other cut is parallel to $CD$. If it contained an adjacent edge, say
$AC$, it would also contain $D$ and would be the face plane $ACD$,
which cannot bisect volume. For opposite edges the two equations are
$\lambda_A=\lambda_B$, $\lambda_C=\lambda_D$. The third depends only
on $s=\lambda_A+\lambda_B$, whose density is $6s(1-s)$; its median is
$s=1/2$. Thus the common point has all barycentric coordinates $1/4$.

\subsection{Three singleton cuts of a regular tetrahedron}
\label{sec:v24:regular-s3}
The three cap-volume formulas below are instances of one change-of-coordinates rule.

\begin{lemma}[A simplex clipping formula]\label{lem:simplex-clipping}
Let $M$ be an invertible $(n+1)\times(n+1)$ matrix with $M^{-1}\ge0$
entrywise, and let $c_j$ be the column sums of $M^{-1}$. In barycentric
coordinates $\lambda$, the region $M\lambda\ge0$ occupies the volume fraction
\begin{equation}\label{eq:simplex-clipping}
 \frac1{|\det M|\prod_{j=1}^{n+1}c_j}
\end{equation}
of the simplex.
\end{lemma}
\begin{proof}
The $c_j$ are positive. Set $y=M\lambda$. Since $M^{-1}\ge0$, the region
becomes $y\ge0$, $\sum c_jy_j=1$. Its vertices are the columns of
$M^{-1}$ divided by their column sums. Their barycentric determinant
is the volume fraction~\eqref{eq:simplex-clipping}.
\end{proof}

For a matrix $B$, write $\Delta_{ij}(B)$ for its principal
two-by-two minor.
\begin{lemma}\label{lem:v24:m-matrix-ratio}
If a real $3\times3$ matrix $B$ has positive diagonal, strictly negative off-diagonal
entries and positive row sums, then its determinant and principal
two-by-two minors are positive, $B^{-1}\ge0$ entrywise, and
\begin{equation}\label{eq:v24:m-ratio}
 \mathcal R(B)=
 \frac{\Delta_{12}(B)\Delta_{13}(B)\Delta_{23}(B)}
      {\det B\,B_{11}B_{22}B_{33}}>1.
\end{equation}
With nonpositive off-diagonal entries the same conclusions hold
with $\mathcal R(B)\ge1$.
\end{lemma}
\begin{proof}
Write $B=D(I-A)$, with $D$ its positive diagonal. Then $A\ge0$,
has zero diagonal and row sums less than one. Its Neumann inverse
is nonnegative, and continuity along $I-tA$ gives positive
determinant. Set $x=A_{12}A_{21}$, $y=A_{13}A_{31}$,
$z=A_{23}A_{32}$, $r=A_{12}A_{23}A_{31}$,
$s=A_{13}A_{32}A_{21}$. The normalized determinant is
$1-x-y-z-r-s>0$, and the normalized product of principal minors
is $(1-x)(1-y)(1-z)$. Their difference
$xy(1-z)+xz+yz+r+s$ is positive in the strict case and
nonnegative in the weak case. Also $x,y,z<1$, proving the
minor signs.
\end{proof}

\begin{proof}[Proof of Theorem~\ref{thm:v24:regular-s3-exclusion}]
Use the vertices $v_1,\ldots,v_4$ of $S_0$ in their displayed
order. They satisfy $v_i\cdot v_j=4\delta_{ij}-1$,
$\sum v_i=0$, $\det[v_1\,v_2\,v_3]=4$, and
$\vol S_0=8/3$. A singleton median with heights $t,-b,-c,-d$
has cap fraction $t^3/((t+b)(t+c)(t+d))=1/2$.
Thus $(1+b/t)(1+c/t)(1+d/t)=2$ and $b+c+d<t$;
its centroid value is positive. Its gradient lies strictly in
the coordinate octant of the isolated vertex: after the
corresponding coordinate sign changes, each component is at
least $(t-b-c-d)/4>0$. Two gradients in that octant cannot
be orthogonal. Apply a tetrahedral symmetry and order the cuts
so that their isolated vertices are the displayed $v_1,v_2,v_3$.

Write $H_i(X)=u_i\cdot X+\eta_i$ with oriented orthonormal
$u_i$ and $\eta_i>0$. Set
\begin{equation}\label{eq:v24:two-matrices}
 \begin{gathered}
 U=(u_i^{\mathsf T})_{i=1}^3,\quad V=[v_1\,v_2\,v_3],
 \quad P=UV,\quad H=P+\eta\mathbf1^{\mathsf T},\\
 p_i=P_{ii},\quad t_i=H_{ii},\quad
 w=P\mathbf1=-Uv_4,\quad h=1+w^{\mathsf T}\eta.
 \end{gathered}
\end{equation}
Both matrices have positive diagonal, negative off-diagonal
entries and positive row sums $w_i,w_i+3\eta_i$, since
$H_i(v_4)=\eta_i-w_i<0$. Lemma~\ref{lem:v24:m-matrix-ratio} applies to both.
It gives $\det P>0$, hence $\det U=1$ and $\det P=4$.
The regular Gram matrix yields
\begin{equation}\label{eq:v24:cofactor-identity}
 \begin{gathered}
 P^{\mathsf T}P=4I-\mathbf1\mathbf1^{\mathsf T},\qquad
 P^{-1}=(P^{\mathsf T}+\mathbf1w^{\mathsf T})/4,\\
 \Delta_{jk}(P)=p_i+w_i,\quad
 \mathcal R(P)=\frac{\prod_i(p_i+w_i)}{4\prod_i p_i}>1,\quad
 \det H=4h.
 \end{gathered}
\end{equation}
Here the last identity uses $\mathbf1^{\mathsf T}P^{-1}=w^{\mathsf T}$.

Let $P_I$ be the fraction where all forms indexed by $I$ are positive,
put $E_{ij}=u_i\cdot(v_i-v_j)>0$ and $\Delta_{ij}=\Delta_{ij}(H)$.
For distinct $i,j,k\in\{1,2,3\}$ the three cell volumes are
\begin{equation}\label{eq:v24:cap-volumes}
 \begin{aligned}
 P_i&=\frac{t_i^3}{E_{ij}E_{ik}(p_i+w_i)},\\
 P_{ij}&=\frac{\Delta_{ij}^3}{16E_{ij}E_{ji}w_kp_k},\\
 P_{123}&=\frac{h^3}{16w_1w_2w_3}
        =\frac{(\det H)^3}{1024w_1w_2w_3}.
 \end{aligned}
\end{equation}
Apply Lemma~\ref{lem:simplex-clipping} with rows given by the selected
forms $H_i$ and the remaining barycentric coordinates. Their leading blocks are principal blocks of $H$, with positive row sums
and nonpositive off-diagonal entries. The Neumann-series argument of
Lemma~\ref{lem:v24:m-matrix-ratio} gives nonnegative inverses; the remaining
columns preserve this property because all other vertex heights are negative.
For one cut the determinant is $t_i$ and the column sums are
$1/t_i$ and $E_{ij}/t_i$ ($j\ne i$), with $E_{i4}=p_i+w_i$.
For two cuts the determinant is $\Delta_{ij}$ and the column sums are
\[
 \frac1{\Delta_{ij}}(E_{ji},E_{ij},B_k,B_4),\qquad
 B_\ell=\Delta_{ij}-E_{ji}H_i(v_\ell)-E_{ij}H_j(v_\ell).
\]
For cyclic $(i,j,k)$, put
$g=E_{ji}u_i+E_{ij}u_j=u_k\times(v_i-v_j)$. Then
$B_\ell=g\cdot(v_i-v_\ell)$, and
\[
 (v_i-v_j)\times(v_i-v_k)=-4v_4,\qquad
 (v_i-v_j)\times(v_i-v_4)=4v_k
\]
give $B_k=4w_k$ and $B_4=4p_k$.
For three cuts the determinant is $4h$ and the column sums are
$(w_1,w_2,w_3,4)/h$: indeed,
$\mathbf1^{\mathsf T}H^{-1}=w^{\mathsf T}/h$, the last column has sum
$1+w^{\mathsf T}(w-\eta)/h=4/h$, and $|w|^2=|v_4|^2=3$.
Substitution into~\eqref{eq:simplex-clipping} gives all three formulas.

Multiplying~\eqref{eq:v24:cap-volumes} and using~\eqref{eq:v24:m-ratio} yields
\begin{equation}\label{eq:v24:multiplicative-identity}
 \frac{P_{123}P_1P_2P_3}{P_{12}P_{13}P_{23}}
 =\frac{4\prod_i p_i}{\prod_i(p_i+w_i)}
  \left(\frac{\det H\prod_i t_i}{\Delta_{12}\Delta_{13}\Delta_{23}}\right)^3
 =\frac1{\mathcal R(P)\mathcal R(H)^3}<1.
\end{equation}
Since $P_i=1/2$ and $P_{ij}=1/4$, the left side is
$8P_{123}<1$, so \eqref{eq:gen:triple-moment} gives $\tau<0$.
\end{proof}

\subsection{Lifting the centered families}
\begin{proposition}\label{prop:liftfamily}
For every $n\ge4$ the regular $n$-simplex admits a nonconstant continuous family of orthogonal partitions.
\end{proposition}
\begin{proof}
For a tetrahedral face $F$ with centroid $g$, put $U=\operatorname{span}(F-g)$.
The remaining vertices are equidistant from those of $F$, hence lie in
$g+U^\perp$. Over a point $y$ of the complementary projection, let $s$
be the sum of the four face barycentric coordinates. Its fiber is
$y+s(F-g)$, so its volume density is proportional to $s^3$. Each centered
triple on $F$ extended along $U^\perp$ has the same sign distribution on
every interior fiber. Its signs are thus independent of the projected point.
Lift a partition of the projection from Theorem~\ref{intro:orthogonal-existence}
(a median for $n=4$). Every mixed pair quarters by independence, and every
unmixed pair quarters by construction. The face triple varies continuously
and is recovered by restricting the three cuts with normals in $U$ to $F$.
\end{proof}

\Needspace{7\baselineskip}
\section{Computational results and further questions}
\label{sec:computations}
\label{sec:computational-results}

The reflection reduction provides a first step toward counting tetrahedral
partitions: its weighted triangle projection admits a complete solution map
(Section~\ref{subsec:weighted-map}). We then consider unrestricted tetrahedral
partitions, obtaining a lower bound of seventeen and a restriction for the
regular shape (Section~\ref{subsec:seventeen-computations}). These results are proved
using exact arithmetic and interval bounds. Their proofs follow in
Section~\ref{subsec:computational-proofs}; the final subsection states the
remaining questions.

\subsection{The weighted triangle map}
\label{subsec:weighted-map}

The reflection reduction in Theorem~\ref{thm:c11-edge} gives the
probability density~\eqref{eq:c11-weight} on
\[
 \Delta(a,H)=\operatorname{conv}\{M=(0,0),C=(1,0),D=(a,H)\},
 \qquad H>0,\qquad K=a^2+H^2.
\]
Its barycentric law is $\operatorname{Dir}(2,1,1)$. Let $N(a,K)$ count its
perpendicular quarterings, for $K>a^2$. The \emph{discriminant} $\Sigma$
is the set of parameters for which a quartering is not simple in the sense
of Section~\ref{gen:definitions}. A vertex incidence means that a cutting
line passes through a vertex. The incidence tip below arises at such a
configuration; the other two tips are ordinary semicubical cusps,
locally equivalent to $(t^2,t^3)$ under smooth changes of coordinates.
Figure~\ref{fig:hz} displays this count in tetrahedral parameters $(h,z)$.

\begin{theorem}[Weighted triangle map, computer-assisted]
\label{thm:weighted-atlas-certified}
The discriminant $\Sigma$ is a simple closed curve with two ordinary
semicubical cusps and an incidence tip at $(1/5,1)$. The complete count is
\[
 N(a,K)=
 \begin{cases}
 3,&\text{in the bounded interior of }\Sigma,\\
 1,&\text{outside its closure},\\
 2,&\text{on }\Sigma\text{ away from the three tips},\\
 1,&\text{at each tip}.
 \end{cases}
\]
Every quartering off $\Sigma$ is simple. In particular, $(a,K)=(1/4,1)$
lies in the open region of three simple quarterings.
\end{theorem}

By Theorem~\ref{thm:c11-edge}, this counts all full tetrahedral partitions
containing the prescribed mirror. Sufficient axial stretch forces every
partition to contain that mirror, giving the exact-count constructions
of Section~\ref{sec:general}.

\subsection{Tetrahedral examples and restrictions}
\label{subsec:seventeen-computations}

Recall from Section~\ref{gen:definitions} that a partition is \emph{regular}
when the three pairwise moment equations have a nonsingular $3\times3$
Jacobian in local rotation coordinates, with each plane kept at its median.
Such a solution is isolated and persists under small shape changes by
Lemma~\ref{gen:continuation}. A small deformation of the regular tetrahedron
gives many such solutions.

\begin{proposition}[Seventeen partitions, computer-assisted]
\label{prop:computational-seventeen}
The tetrahedron
\[
 T_* =\conv\{(-1,0,0),(1,0,0),(0,\sqrt3,0),
          (0,1/\sqrt3,\sqrt{13317/5000})\}
\]
has at least seventeen distinct regular orthogonal partitions, including
three full partitions. Every sufficiently small perturbation of its vertices
still has at least seventeen regular orthogonal partitions.
\end{proposition}

For the regular shape, the following restriction supplements the centered
classification of Section~\ref{sec:tetrahedra}.

\begin{proposition}[One central cut, computer-assisted]
\label{prop:computational-central-cut}
Every orthogonal partition of the regular tetrahedron with a plane through
its centroid is centered. Any full octasection beyond the three centered
ones has three noncentral planes.
\end{proposition}

\subsection{Proofs and exact computations}
\label{subsec:computational-proofs}

The computational supplement contains the rational enclosures and full
algebraic reductions. Existence and uniqueness of local roots follow from a common contraction test.
For a square $C^1$ system $F$, a box $X=c+[-\rho,\rho]^d$ with $\rho>0$, and a fixed matrix $R$, the test
\begin{equation}\label{eq:comp:box-test}
 \|RF(c)\|_\infty+\rho\sup_X\|I-RDF\|_\infty<\rho
\end{equation}
guarantees a unique regular zero in $X$. Indeed, $x\mapsto x-RF(x)$ contracts
$X$ into itself, and the derivative bound makes $RDF$ invertible by the
Neumann series. Thus $R$ is invertible and its fixed point is a regular zero.

\begin{proof}[Proof of Theorem~\ref{thm:weighted-atlas-certified} (outline)]
All affine quarterings form three compact arcs, indexed by the side bearing
two chord endpoints, with shared endpoints identified. Monotonicity of cap
mass gives one branch per arrangement; an overlap bound excludes support
on only two sides. With arc parameter $s$, the orthogonality equation
$F(s;a,K)=0$ is linear in $(a,K)$. Its envelope is given by
$F=\partial_sF=0$, retaining the undivided
equations at vertex incidences.

Exact interval bounds establish one ordinary cusp on each of two symmetric
envelopes, monotonicity of the third, and separation apart from the prescribed
joins. Hence their union is a Jordan curve. The boundary-annulus argument of
Lemma~\ref{gen:planar-seeds} gives count one outside, and the same lemma
supplies a region of count three. Continuation therefore gives count three
throughout the interior. At a smooth fold two roots merge; at an ordinary
cusp three merge. A separate incidence calculation gives one root at the
remaining tip. These checks determine all boundary counts.

At $(a,K)=(1/4,1)$, normalize the isosceles base $CD$ to length one. Its
squared altitude satisfies
\[
 \frac38<k=\frac{1+a}{4(1-a)}=\frac5{12}
       <\frac1{\sqrt2}-\frac14.
\]
The monotonicity bounds give one root on each lateral arc. The symmetry
axis and $\lambda_M=1/2$ give the third, since
$\lambda_M\sim\operatorname{Beta}(2,2)$. The count table places the point
in the three-root interior.
\end{proof}

\begin{proof}[Proof of Proposition~\ref{prop:computational-seventeen}]
Four configurations are located and verified in normalized vertex-height coordinates.
For the three without vertex incidences, the test~\eqref{eq:comp:box-test} applies to nine median,
quarter-mass and orthogonality equations. The reflection configuration is
verified in its symmetry slice; section moments establish nonsingularity of
its unrestricted Jacobian. Thus all four configurations are regular
by Lemma~\ref{gen:cut-flux}; the clipping inequalities are also verified.
Under the six base isometries, the reflection, two singleton, and two--two
configurations have stabilizer orders $2,1,3,1$. Their orbits therefore
give $3+6+2+6=17$ roots. Interval separation
checks distinguish every pair under every matching of the three unoriented
planes. The first orbit contains an edge-containing cut and is full by
Theorem~\ref{thm:c11-edge}. Regularity and Lemma~\ref{gen:continuation} give
persistence under arbitrary small shape perturbations.
\end{proof}

\begin{proof}[Proof of Proposition~\ref{prop:computational-central-cut} (outline)]
Lemma~\ref{c11:median} and symmetry normalize the central cut to $x+ty=0$,
$0\le t\le1$. One remaining median can then be written
$-tx+y+vz=m$, with $0\le v\le\sqrt{1+t^2}$. The cases $t=0$ and $v=0$
are centered by the same lemma. At $t=1$, the reflection projection has
$(a,K)=(1/3,1)$; the exact count is one, attained by a centered partition.
For $0<t<1$ and $v>0$, the actual median has two possible vertex-sign
patterns. In each, eliminating $m$ from the median and quartering equations
gives a polynomial in $(t,v)$. Rational Bernstein coefficients prove its
strict positivity on the entire admissible region, including the common
incidence boundary. Thus no further quartering is possible in this sector.
\end{proof}

The corresponding supplement files are grouped in \texttt{research/v10/}
(weighted map), \texttt{research/v18/} (seventeen roots), and
\texttt{research/v32/} (central cut).

\subsection{Further questions}
\label{subsec:computational-questions}

\paragraph{Full partitions.}
Do the three centered octasections exhaust all full partitions of the regular
tetrahedron? Proposition~\ref{prop:computational-central-cut} restricts an
additional partition to three noncentral cuts, while
Theorem~\ref{thm:v24:regular-s3-exclusion} excludes three strict singleton cuts.
The other cut types and their incidences remain to be treated.
In higher dimensions, the first case beyond the binary criterion is the regular
$14$-simplex, since $15=(1111)_2$: does it admit a full partition?
Generic existence of full partitions is open even for tetrahedra. A dimension
argument must control the strata on which the additional fullness equation
could vanish identically.

\paragraph{Counts and continuous families.}
What is the sharp upper bound for the finite generic tetrahedral count,
and how does it change between regions of shape space? Can a fixed tetrahedron
have a nonconstant family outside its centered locus? Generic finiteness and
the centered classification leave both questions open.

\paragraph{Further computations.}
The weighted map suggests analogous studies for other Dirichlet parameters
and in higher dimensions. Complete counts require a description of all
admissible affine configurations as well as their local singularities.
Effective stretch bounds in arbitrary dimension would also make the
exact-count constructions explicit.

\FloatBarrier

\paragraph{Acknowledgments.}
This work was carried out during my sabbatical at EPFL in Lausanne.
I am grateful to Maryna Viazovska for inviting me and for the opportunity
to work at EPFL. The Bernoulli Center for Fundamental Studies, on whose
Executive Committee she serves, was one of the reasons I was drawn to
Jacob Bernoulli's quadrisection problem.

\paragraph{Use of AI.}
The author used tools from OpenAI for assistance with the computations,
the exposition of the proofs, and the preparation of the manuscript.
The author is responsible for the arguments and the final text.

\bigskip
\noindent\textsc{Oleg R. Musin}\\
School of Mathematical and Statistical Sciences\\
The University of Texas Rio Grande Valley\\
One West University Boulevard, Brownsville, TX 78520, USA\\
\textit{E-mail:} \href{mailto:oleg.musin@utrgv.edu}{oleg.musin@utrgv.edu}
\end{document}